\documentclass[10pt,a4paper]{amsart}
\usepackage{graphicx,latexsym,amssymb,diagrams}

\oddsidemargin=\evensidemargin
\catcode`\@=11
\long\def\@makefntext#1{
\protect\noindent \hbox to 3.2pt {\hskip-.9pt  
$^{{\eightrm\@thefnmark}}$\hfil}#1\hfill}		

\def\@makefnmark{\hbox to 0pt{$^{\@thefnmark}$\hss}}	
	
\def\ps@myheadings{\let\@mkboth\@gobbletwo		
\def\@oddhead{\hbox{}
\rightmark\hfil\eightrm\thepage}   
\def\@oddfoot{}\def\@evenhead{\eightrm\thepage\hfil
\leftmark\hbox{}}\def\@evenfoot{}
\def\sectionmark##1{}\def\subsectionmark##1{}}

\catcode`@=11		      		     
\def\ps@plain{\let\@mkboth\@gobbletwo
     \def\@oddhead{}\def\@oddfoot{\eightrm\hfil\thepage
     \hfil}\def\@evenhead{}\let\@evenfoot\@oddfoot}

\newcounter{sectionc}\newcounter{subsectionc}\newcounter{subsubsectionc}
\renewcommand{\section}[1] {\vspace{12pt}\addtocounter{sectionc}{1} 
\setcounter{subsectionc}{0}\setcounter{subsubsectionc}{0}\noindent 
	{\tenbf\thesectionc. #1}\par\vspace{5pt}}
\renewcommand{\subsection}[1] {\vspace{12pt}\addtocounter{subsectionc}{1} 
	\setcounter{subsubsectionc}{0}\noindent 
	{\bf\thesectionc.\thesubsectionc. 
	{\kern1pt \bfit #1}}\par\vspace{5pt}}
\renewcommand{\subsubsection}[1] {\vspace{12pt}
	\addtocounter{subsubsectionc}{1}
	\noindent
	{\tenrm\thesectionc.\thesubsectionc.\thesubsubsectionc.	{\kern1pt 
	\it #1}}\par\vspace{5pt}}

\newcommand{\textlineskip}{\baselineskip=13pt}
\newcommand{\smalllineskip}{\baselineskip=10pt}

\newcommand{\copyrightheading}[1]
	{\vspace*{-2.5cm}\smalllineskip{\flushleft
	{\footnotesize Journal of Knot Theory and Its Ramifications #1}\\
   	{\footnotesize \copyright\kern2pt World Scientific 
        Publishing Company}\\
         }}

\def\abstracts#1#2#3#4{{
	\centering{\begin{minipage}{4.5in}\footnotesize\baselineskip=10pt
	\centerline{ABSTRACT} 
	\parindent=15pt #1\par 
	\parindent=15pt #2\par
	\parindent=15pt #3\par
	\parindent=15pt #4\par
	\end{minipage}}\par}} 

\def\keywords#1{{ 
	\centering{\begin{minipage}{4.5in}\footnotesize\baselineskip=10pt
	{\footnotesize\it Keywords}\/: #1
	\end{minipage}}\par}}

\def\dates#1{{ 
	\centering{\begin{minipage}{4.5in}\footnotesize\baselineskip=10pt
	{\footnotesize\it Date}\/: #1
	\end{minipage}}\par}}

\renewenvironment{thebibliography}[1]
	{\frenchspacing
	 \ninerm\baselineskip=11pt
	 \begin{list}{[\arabic{enumi}]}
	{\usecounter{enumi}\setlength{\parsep}{0pt}
	 \setlength{\leftmargin 13.7pt}{\rightmargin 0pt} 
	 \setlength{\leftmargin 19pt}{\rightmargin 0pt}   
	 \setlength{\leftmargin 24pt}{\rightmargin 0pt}   
	 \setlength{\itemsep}{0pt} \settowidth
	{\labelwidth}{[#1]}\sloppy}}{\end{list}}

\newcounter{itemlistc}
\newcounter{romanlistc}
\newcounter{alphlistc}
\newcounter{arabiclistc}

\newcommand{\fcaption}[1]{
        \refstepcounter{figure}
        \setbox\@tempboxa = \hbox{\footnotesize Fig.~\thefigure. #1}
        \ifdim \wd\@tempboxa > 5in
           {\begin{center}
        \parbox{5in}{\footnotesize\smalllineskip Fig.~\thefigure. #1}
            \end{center}}
        \else
             {\begin{center}
             {\footnotesize Fig.~\thefigure. #1}
              \end{center}}
        \fi}
\def\pmb#1{\setbox0=\hbox{#1}
	\kern-.025em\copy0\kern-\wd0
	\kern.05em\copy0\kern-\wd0
	\kern-.025em\raise.0433em\box0}

\def\fpage#1{\begingroup
\voffset=.3in
\thispagestyle{empty}
\endgroup}

\font\tenrm=cmr10
 
\font\tenbf=cmbx10
\font\bfit=cmbxti10 at 10pt
\font\ninerm=cmr9
\font\nineit=cmti9

\font\eightrm=cmr8

\theoremstyle{plain} 
\newtheorem{theorem}{Theorem}[sectionc]
\newtheorem{lemma}[theorem]{Lemma}
\newtheorem{proposition}[theorem]{Proposition}
\newtheorem{corollary}[theorem]{Corollary}
\theoremstyle{definition}

\newtheorem{remark}[theorem]{Remark}

\newtheorem{example}[theorem]{Example} 
\newtheorem{conventions}[theorem]{Convention} 
\numberwithin{equation}{sectionc}
    	{\setcounter{itemlistc}{0}		
	 \begin{list}{$\bullet$}		
	{\usecounter{itemlistc}			
	 \leftmargin10pt	       
	 \setlength{\parsep}{0pt}
	 \setlength{\itemsep}{0pt}     
	}}{\end{list}}

	{\setcounter{romanlistc}{0}		
	 \begin{list}{$($\roman{romanlistc}$)$}	
	{\usecounter{romanlistc}		
	 \leftmargin18pt 
	 \setlength{\parsep}{0pt}
	 \setlength{\itemsep}{0pt}	
	 \settowidth{\labelwidth}{#1}                          
	}}{\end{list}}

	{\setcounter{enumii}{0}			
	 \begin{list}{$($\alph{enumii}$)$}	
	{\usecounter{enumii}			
	 \leftmargin18pt		
	 \setlength{\parsep}{0pt}
	 \setlength{\itemsep}{0pt}	
	 \settowidth{\labelwidth}{#1}                          
	}}{\end{list}}

\def\qed{\hbox{${\vcenter{\vbox{			
   \hrule height 0.4pt\hbox{\vrule width 0.4pt height 6pt
   \kern5pt\vrule width 0.4pt}\hrule height 0.4pt}}}$}}

\newcommand{\bd}{\begin{description}}   
\newcommand{\ed}{\end{description}} 
\newcommand{\ba}{\begin{array}}      \newcommand{\ea}{\end{array}} 
\newcommand{\bc}{\begin{center}}     \newcommand{\ec}{\end{center}} 
\newcommand{\be}{\begin{enumerate}}  \newcommand{\ee}{\end{enumerate}} 
\newcommand{\beq}{\begin{eqnarray}}  \newcommand{\eeq}{\end{eqnarray}} 
\newcommand{\beQ}{\begin{eqnarray*}} \newcommand{\eeQ}{\end{eqnarray*}} 
\newcommand{\bi}{\begin{itemize}}    \newcommand{\ei}{\end{itemize}}

\newcommand{\ov}{\overline} 
\newcommand{\we}{\wedge} 
\newcommand{\Y}{\mathsf{Y}} 
\newarrow{Noto} ----- 
\newarrow{To} ----> 
\newarrow{Mapsto} |---> 
\newarrow{Into} >---> 
\newarrow{Onto} ----{>>} 
\newarrow{Iso} >---{>>} 
\newarrow{Inclus} C---> 
\newarrow{Equal} ===== 
\newarrow{Dash} {dash}{dash}{dash}{dash}{dash} 
\newarrow{Dashto} {dash}{dash}{dash}{dash}{>} 
\begin{document}
\setlength{\textheight}{7.7truein}  


\normalsize\textlineskip
\setcounter{page}{1}


\vspace*{0.88truein}

\fpage{1}
\centerline{\bf CHARACTERIZATION OF $Y_2$-EQUIVALENCE}
\baselineskip=13pt
\centerline{\bf FOR HOMOLOGY CYLINDERS}
\vspace*{0.37truein}
\centerline{\footnotesize GW\'ENA\"EL MASSUYEAU and JEAN-BAPTISTE MEILHAN}
\baselineskip=12pt
\centerline{\footnotesize\it Laboratoire Jean Leray, UMR 6629 CNRS/Universit\'e de Nantes}
\baselineskip=10pt
\centerline{\footnotesize\it 2 rue de la Houssini\`ere, BP 92208, 	
             44322 Nantes Cedex 03, France}
\baselineskip=12pt
\centerline{\footnotesize\it massuyea@math.univ-nantes.fr, meilhan@math.univ-nantes.fr}

\vspace{1cm}

\abstracts{ For $\Sigma$ a compact connected oriented surface, we consider  
homology cylinders  over $\Sigma$: these are homology  
cobordisms with an extra homological triviality  
condition. When considered up to $Y_2$-equivalence,  
which is a surgery equivalence relation arising from the  
Goussarov-Habiro theory, homology cylinders form an Abelian group.\\
In this paper, when $\Sigma$ has one or zero boundary
component, we define a surgery map from a certain space of 
graphs to this group. This map is shown to be an isomorphism, 
with inverse given by some  extensions of the first Johnson 
homomorphism and Birman-Craggs homomorphisms.}{}{}{}

\vspace*{10pt}
\dates{March 18, 2002 and, in revised form, September  5, 2002}


%
%
%
\section{Introduction}	
\label{sec:intro}

\subsection{Homology cylinders} 
Homology cylinders are important objects  
in the theory of finite type invariants of Goussarov-Habiro: 
they have thus appeared in both \cite{H} and \cite{G}.  
Let us recall the definition of these objects. 
 
Let $\Sigma$ be  a compact connected oriented surface. 
A \emph{homology cobordism} over $\Sigma$ is a triple $(M,i^{+},i^{-})$ 
where $M$ is a compact oriented $3$-manifold and  
$i^{\pm}: \Sigma \rTo M$ are oriented embeddings with images $\Sigma^\pm$, such that:
\begin{itemize}
\item[(i)] $i^\pm$ are homology isomorphisms;
\item[(ii)] $\partial M=\Sigma^+\cup \left(-\Sigma^-\right)$ and  
$\Sigma^+\cap\left(-\Sigma^-\right)=\pm\partial\Sigma^\pm$;
\item[(iii)] $i^+|_{\partial \Sigma}=i^-|_{\partial \Sigma}$.
\end{itemize} 
Homology cobordisms are considered up to orientation-preserving diffeomorphisms. 
When $(i^-)_*^{-1}\circ (i^+)_*: H_1(\Sigma;\mathbf{Z})\rTo H_1(\Sigma;\mathbf{Z})$
is the identity, $M$ is said to be a \emph{homology cylinder}.  
The set of homology cobordisms  
is denoted here by $\mathcal{C}(\Sigma)$, 
and $\mathcal{HC}(\Sigma)$ denotes the subset of homology cylinders.  
If $M=(M,i^+,i^-)$ and $N=(N,j^+,j^-)$ are 
homology cobordisms, we can define their \emph{stacking product} by
\begin{displaymath}
M\cdot N := (M\cup_{i^-\circ (j^+)^{-1}}N,i^+,j^-). 
\end{displaymath}
This product induces a monoid structure on $\mathcal{C}(\Sigma)$, 
with $\mathcal{HC}(\Sigma)$ a submonoid.  
The unit element is $1_\Sigma:=\left(\Sigma \times I,Id,Id\right)$, 
where $I$ is the unit interval $[0,1]$ and where a collar of $\Sigma^{\pm}$  
is stretched along $\partial \Sigma \times I$ so that the second defining condition
for homology cobordisms is satisfied. 
 
Habiro in \cite[\S8.5]{H} outlined how homology cylinders  
can serve as a powerful tool  
in studying the mapping class groups of 
surfaces (see \cite{GL}, \cite{Ha}, \cite{L}). 
The connection lies on the homomorphism of monoids 
$$
\mathcal{T}(\Sigma)  \rTo^{C}  \mathcal{HC}(\Sigma)  
$$ 
sending each $h$ in the Torelli group of $\Sigma$ to the mapping 
cylinder $C_h=(\Sigma\times I, Id, h)$ (with, as above, a collar of $\Sigma^{\pm}$ 
stretched along $\partial \Sigma \times I$).\\
In the sequel, \emph{we restrict ourselves} to the following two cases:
\begin{itemize}  
\item[(i)] $\Sigma=\Sigma_g$ is the standard closed 
oriented surface of genus $g\geq 0$, which here is referred 
to as the \emph{closed case}; 
\item[(ii)] $\Sigma=\Sigma_{g,1}$ is the standard   
compact oriented surface of genus $g\geq 0$ with  
one boundary component, which here is referred to as the \emph{boundary case}.  
\end{itemize}
The usual notations $\mathcal{T}_{g,1}=\mathcal{T}(\Sigma_{g,1})$ 
and $\mathcal{T}_{g}=\mathcal{T}(\Sigma_{g})$ 
for the Torelli groups will be used. Also denote by  
$H$ the first homology group of $\Sigma$ with  
integer coefficients, by $\bullet$ the intersection form on $H$ 
and by $(x_i,y_i)_{i=1}^g$ a symplectic basis for $(H,\bullet)$.

\subsection{$Y_k$-\textbf{\textit{equivalence}}}
The theory of finite type invariants  
of Goussarov-Habiro has come equipped with a  
topological calculus toolbox:  
this was called \emph{calculus of claspers} in  
\cite{H} or alternatively \emph{clovers} in \cite{GGP}.  
We will assume a certain familiarity of the reader with these techniques.\\
In particular, let us recall that, for $k \geq 1$ an  
integer, the \emph{$Y_k$-equivalence}\footnote{This equivalence  
relation is called $(k-1)$-equivalence in \cite{G},  
and $A_k$-equivalence in \cite{H}.}\,\, is the equivalence 
relation generated by surgery on connected clovers of degree $k$. 
Following Habiro in \cite{H}, we can then define  
a descending filtration of monoids 
\begin{displaymath} 
\mathcal{C}(\Sigma) \supset \mathcal{C}_1(\Sigma) \supset  
\mathcal{C}_2(\Sigma) \supset \cdots \supset  
\mathcal{C}_k(\Sigma) \supset \cdots 
\end{displaymath} 
where $\mathcal{C}_k(\Sigma)$ is the submonoid  
consisting of the homology cobordisms which are  
$Y_k$-equivalent to the trivial cobordism $1_\Sigma$.
Note the following fact,
a proof of which has been inserted in \S 4.

\begin{proposition}
\label{prop:homology_cylinders}
If $\Sigma=\Sigma_g$ or $\Sigma_{g,1}$, then
$\mathcal{HC}(\Sigma)=
\mathcal{C}_1(\Sigma)$.
\end{proposition}

\noindent
As mentioned by Habiro, we can show from the calculus of 
clovers that for every $k\geq 1$, the quotient monoid  
\begin{displaymath} 
\overline{\mathcal{C}}_k(\Sigma):=\mathcal{C}_k(\Sigma)/Y_{k+1} 
\end{displaymath} 
is an Abelian group. In particular, $\overline{\mathcal{C}}_1(\Sigma)$  
is the Abelian group of homology cylinders over $\Sigma$  up to $Y_2$-equivalence.
This group is the subject of the present paper. 

For $k\geq 2$, Habiro gives a combinatorial  
upper bound for the Abelian group
$\overline{\mathcal{C}}_k(\Sigma)$. Precisely,  
he defines $\mathcal{A}_{k}(H)$  
to be the Abelian group (finitely) generated by unitrivalent graphs  
of internal degree $k$, with cyclic orientation at each trivalent vertex 
and whose univalent vertices are labelled by elements of $H$  
and are totally ordered. These graphs are considered modulo  
the well-known AS, IHX, multilinearity 
relations, and up to some ``STU-like relations''  
dealing with the order of the univalent 
vertices. In the closed case, some relations of a symplectic type  
can be added. Then, there is a surjective \emph{surgery map} 
$$
\mathcal{A}_k(H)  \rOnto^{\psi_{k}}   \overline{\mathcal{C}}_k(\Sigma) 
$$
sending each graph $G$ to $(1_\Sigma)_{\tilde{G}}$,  
where $\tilde{G}$ is a clover in  
the manifold $1_\Sigma$ with $G$ as associated abstract graph,  
whose leaves are stacked from the upper surface $\Sigma \times 1$ 
according to the total order, framed along this surface and   
embedded  according to the labels  
of the corresponding univalent vertices.  
The fact that $\psi_{k}$  is well-defined also 
follows from the calculus of clovers. 
 
As for the case $k=1$, Habiro does not define  
any space of graphs but announces the following isomorphisms
\begin{equation}
\label{eq:iso_Habiro}
\left\{\begin{array}{l}
\overline{\mathcal{C}}_1(\Sigma_{g,1})\simeq
\Lambda^3H \oplus \Lambda^2H_{(2)} \oplus H_{(2)} \oplus \mathbf{Z}_2 \\
\overline{\mathcal{C}}_1(\Sigma_{g})\simeq 
\Lambda^3H/(\omega \wedge H ) \oplus \Lambda^2H_{(2)}/\omega_{(2)}  
\oplus H_{(2)} \oplus \mathbf{Z}_2
\end{array} \right.
\end{equation}
where $H_{(2)}= H\otimes \mathbf{Z}_2$ and where
$$\omega=\sum_{i=1}^g x_i\wedge y_i \in \Lambda^2H$$ 
is the \emph{symplectic element}. This fact has been used afterwards in \cite{L}.\\ 
\emph{The goal of this paper is to prove these isomorphisms},  
in a diagrammatic way, by again defining a surgery map
$$
\mathcal{A}_1(P)  \rTo^{\psi_1}  \overline{\mathcal{C}}_1(\Sigma).  
$$ 
The space of graphs $\mathcal{A}_1(P)$ and the map $\psi_1$
appear to be meaningfully different from $\mathcal{A}_k(H)$ and $\psi_k$ for $k>1$, 
making thus the case $k=1$ exceptional. Indeed, their definition will involve 
both the homology group $H$ and $Spin\left(\Sigma\right)$, the set of \emph{spin structures} on $\Sigma$.  

\subsection{The Abelianized Torelli group}
\label{subsec:Torelli}
We denote by $\Omega_g$ the set of \emph{quadratic forms} with 
$\bullet: H_{(2)}\times H_{(2)} \rTo \mathbf{Z}_2$ 
as associated bilinear form, namely  
\begin{displaymath}
\Omega_g=\left\{H_{(2)} \rTo^q  \mathbf{Z}_2 \ : \
\forall x,y \in H_{(2)},\ q(x+y)-q(x)-q(y)=x\bullet y \right\}.
\end{displaymath}
Note that $\Omega_g$ is an affine space 
over $H_{(2)}$, with action given by 
\begin{displaymath} 
\forall q\in \Omega_g, \forall x\in H_{(2)}, \quad 
x\cdot q:= q + x\bullet(-). 
\end{displaymath} 
Thus, among the maps $\Omega_g \rTo  \mathbf{Z}_2$, 
there are the affine functions, and more generally there are 
the \emph{Boolean polynomials} which are defined to be sums of  
products of affine ones (see \cite[\S4]{JBC}). 
These polynomials form a  
$\mathbf{Z}_2$-algebra denoted by $B_g$, which is filtered by the 
\emph{degree} (defined in the obvious way): 
\begin{displaymath} 
B^{(0)}_g \subset B^{(1)}_g \subset \cdots \subset B_g. 
\end{displaymath} 
For instance, $B^{(1)}_g$ is the space of affine functions on $\Omega_g$;
the constant function $\overline{1}: \Omega_g \rTo\mathbf{Z}_{2}$ sending each $q$ to $1$ 
and, for $h\in H$, the function $\overline{h}$ sending each $q$ to $q(h)$ 
are affine functions. Note the following identity:
\begin{equation}
\label{eq:quadratic_identity}
\forall h_1, h_2 \in H, \quad
\overline{h_1+h_2}=\overline{h_1}+\overline{h_2}
+ (h_1 \bullet h_2) \cdot \overline{1} \in B_g^{(1)}.
\end{equation}
Another example of Boolean polynomial is the 
quadratic Boolean function
\begin{displaymath}
\alpha=\sum_{i=1}^g \overline{x_i}\cdot\overline{y_i},
\end{displaymath}
which is known as the \emph{Arf invariant}.
For any basis $(e_i)_{i=1}^{2g}$ for $H$, there is an isomorphism of algebras:
\begin{equation}
\label{eq:iso_Boolean}
B_g \simeq \frac{\mathbf{Z}_2[t_1,\dots,t_{2g}]}
{t_i^2=t_i}
\end{equation}
sending $\overline{1}$ to $1$ and $\overline{e_i}$ to $t_i$.

Recall now from \cite{JBC}, that the many \emph{Birman-Craggs homomorphisms}
can be summed up into a single homomorphism
\begin{diagram}
& \mathcal{T}_{g,1} & \rTo^\beta & B^{(3)}_g &  
\textrm{ or }& \mathcal{T}_g & \rTo^\beta & \frac{B^{(3)}_g}{\alpha\cdot B^{(1)}_g},
\end{diagram}
according to whether one is considering the boundary case or the closed case.
Recall also from \cite{Jab} that the \emph{first Johnson homomorphism} is a homomorphism
\begin{diagram}
 & \mathcal{T}_{g,1} & \rTo^{\eta_1} & \Lambda^3 H & 
\textrm{ or } & \mathcal{T}_{g} & \rTo^{\eta_1} & \frac{\Lambda^3 H}{\omega\wedge H}.
\end{diagram}
Form the following pull-back:
\begin{diagram}
\Lambda^{3} H \times_{\Lambda^{3} H_{(2)}}\SEpbk B^{(3)}_g   
& \rTo &  B^{(3)}_g \\ 
 \dTo &  & \dTo_q \\ 
\Lambda^{3} H  & \rTo_{-\otimes\mathbf{Z}_2} & \Lambda^{3} H_{(2)},  
\end{diagram} 
where the map $q$ is the canonical projection
$ B^{(3)}_g \rTo B^{(3)}_g/ B^{(2)}_g$ followed by the isomorphism
$B^{(3)}_g/ B^{(2)}_g\simeq \Lambda^{3} H_{(2)}$ which identifies
the cubic polynomial $\overline{h_1} \overline{h_2} \overline{h_3}$
with $h_1\wedge h_2\wedge h_3$ (this is well-defined because
of (\ref{eq:quadratic_identity}) and (\ref{eq:iso_Boolean})).\\
We denote by S the subgroup of this pull-back corresponding to 
$\omega\wedge H \subset \Lambda^{3} H$ and 
$\alpha\cdot B^{(1)}_g \subset B^{(3)}_g$.
Johnson has shown  in \cite{JFinal} that, under the assumption $g\geq 3$, 
the homomorphisms $\eta_1$ and $\beta$ induce isomorphisms
\begin{diagram}
&\frac{\mathcal{T}_{g,1}}{\mathcal{T}_{g,1}'}& 
\rTo^{(\eta_1,\beta)}_{\simeq} & 
\Lambda^{3} H \times_{\Lambda^{3} H_{(2)}} B^{(3)}_g &
\quad \textrm{ and }\quad &
\frac{\mathcal{T}_{g}}{\mathcal{T}_{g}'}& 
\rTo^{(\eta_1,\beta)}_{\simeq} & 
\frac{\Lambda^{3} H \times_{\Lambda^{3} H_{(2)}} B^{(3)}_g}{S}. 
\end{diagram}
\begin{remark}
\label{rem:non-canonical}
Note that, because of (\ref{eq:iso_Boolean}),
the codomains of these maps are respectively non-canonically
isomorphic to
$\Lambda^3H \oplus \Lambda^2H_{(2)} \oplus H_{(2)} \oplus \mathbf{Z}_2$ 
and $\Lambda^3H/(\omega \wedge H ) \oplus \Lambda^2H_{(2)}/\omega_{(2)}  
\oplus H_{(2)} \oplus \mathbf{Z}_2$.
\end{remark}

\subsection{Statement of the results} 
In \S 2, we will construct the space of graphs $\mathcal{A}_1(P)$ and
the surgery map $\psi_1: \mathcal{A}_1(P) \rTo \overline{\mathcal{C}}_1\left(\Sigma\right)$.
Spin structures play a prominent role in their definitions.\\
Observe that, $\overline{\mathcal{C}}_1(\Sigma)$ being an Abelian group, the 
mapping cylinder construction induces a group homomorphism 
\begin{diagram} 
\frac{\mathcal{T}(\Sigma)}{\mathcal{T}(\Sigma)'} & 
\rTo^C & \overline{\mathcal{C}}_1(\Sigma). 
\end{diagram}  
As pointed out by Garoufalidis and Levine in \cite{GL} and \cite{L},
Johnson homomorphisms and Birman-Craggs homomorphisms factor
through  $C: \mathcal{T}(\Sigma) \rTo \mathcal{HC}(\Sigma)$.
These extensions will be detailed in \S 3.\\
Next, we will specify in \S 4 an isomorphism 
$\rho: \mathcal{A}_1(P) \rTo \Lambda^{3} H \times_{\Lambda^{3} H_{(2)}} B^{(3)}_g$
and the following two theorems will be proved from the previous material. 
\begin{theorem} 
\label{th:boundary} 
In the boundary case, the diagram  
\begin{diagram}  
\mathcal{A}_1(P) & \rTo^{\psi_1} &
\overline{\mathcal{C}}_1(\Sigma_{g,1}) & \lTo^{C} &  
\frac{\mathcal{T}_{g,1}}{\mathcal{T}_{g,1}'}\\
& \rdTo<{\rho}  & \dTo>{(\eta_{1},\beta)} & \ldTo>{(\eta_{1},\beta)} &\\  
& & \Lambda^{3} H \times_{\Lambda^{3} H_{(2)}} B^{(3)}_g & &  
\end{diagram}
commutes and all of its arrows are \emph{isomorphisms}, 
except for the two maps starting from $\mathcal{T}_{g,1}/\mathcal{T}'_{g,1}$ when $g<3$.
\end{theorem} 
\begin{theorem} 
\label{th:closed} 
In the closed case, the diagram  
\begin{diagram}  
\frac{\mathcal{A}_1(P)}{\rho^{-1}(\textrm{S})}\quad & \rTo^{\psi_1}
&\overline{\mathcal{C}}_1(\Sigma_{g})&  \lTo^{C} & 
\frac{\mathcal{T}_{g}}{\mathcal{T}_g'}\\
& \rdTo<{\rho}  & \dTo>{(\eta_{1},\beta)} & \ldTo>{(\eta_{1},\beta)} &\\ 
& &  \quad \frac{\Lambda^{3} H \times_{\Lambda^{3} H_{(2)}}  B^{(3)}_g}
{\textrm{S}}\quad
& &  
\end{diagram} 
commutes and all of its arrows are \emph{isomorphisms},
except for the two maps starting from $\mathcal{T}_{g}/\mathcal{T}'_{g}$ when $g<3$.
\end{theorem}
\noindent
Note that Theorem \ref{th:boundary} and Theorem \ref{th:closed}
together with Remark \ref{rem:non-canonical}, give Habiro's
isomorphisms (\ref{eq:iso_Habiro}), which are non-canonical.
Also, we will easily deduce the following. 
\begin{corollary} \label{cor:jbc}
For $\Sigma=\Sigma_{g,1}$ or $\Sigma_g$, let
$M$ and $M'$ be two homology cylinders over $\Sigma$.
Then, the following assertions are equivalent:
\begin{enumerate}
\item[(a)] $M$ and $M'$ are  $Y_2$-equivalent;
\item[(b)] $M$ and $M'$ are not distinguished by degree $1$ 
Goussarov-Habiro finite type invariants;
\item[(c)] $M$ and $M'$ are not distinguished by the first Johnson
homomorphism nor Birman-Craggs homomorphisms.
\end{enumerate}
\end{corollary} 
\noindent
Finally, if an embedding $\Sigma_{g,1}\rInclus \Sigma_g$ is fixed,
there is an obvious ``filling-up'' map 
$\mathcal{C}_{1}\left(\Sigma_{g,1}\right) \rTo
\mathcal{C}_{1}\left(\Sigma_{g}\right)$, through 
which the commutative diagrams  
of Theorem \ref{th:boundary} and Theorem \ref{th:closed} are compatible. 
The reader is referred to \S 4 for a precise statement.

\section{Definition of the Surgery Map $\psi_1$} 
\label{sec:definitions}
In this section, we define the space of graphs $\mathcal{A}_1(P)$ and 
the surgery map $\psi_1$ announced in the introduction.

\subsection{Special Abelian groups and the $\mathcal{A}_1$ functor} 
\label{subsec:A1} 
Let us denote by $\mathcal{A}b$ the category of Abelian groups.  
An \emph{Abelian group with special element} is a pair $(G,s)$ 
where $G$ is an Abelian group and $s\in G$ is of order at most $2$.   
We denote by $\mathcal{A}b_s$ the category of special Abelian 
groups whose morphisms are group homomorphisms preserving the special 
elements. We now define a functor 
$$
\mathcal{A}b_s  \rTo^{\mathcal{A}_1}  \mathcal{A}b 
$$
in the following way. For $(G,s)$ an object  
in $\mathcal{A}b_{s}$, $\mathcal{A}_1(G,s)$ is the free Abelian group  
generated by Y-shaped unitrivalent graphs, whose trivalent vertex is  
equipped with a cyclic order on the incident edges  
and whose univalent vertices are labelled by $G$, subject to some relations.  
The notation  
 \[ \Y[z_1,z_2,z_3] \] 
will stand for the $Y$-shaped graph  
whose univalent vertices are colored by $z_1$,  
$z_2$ and $z_3\in G$ in accordance with the cyclic order,  
so that our notation is invariant  
under cyclic permutation of the $z_i$'s. 
The relations are the following ones:\\[0.5cm] 
\begin{tabular}{rcl}  
\textbf{Antisymetry (AS)} & : &
$\Y[z_1,z_2,z_3] =  -\Y[z_2,z_1,z_3],$\\[0.3cm] 
\textbf{Multilinearity of colors} & : &  
$\Y[z_0+z_1,z_2,z_3] = \Y[z_0,z_2,z_3] +  
\Y[z_1,z_2,z_3],$ \\[0.3cm] 
\textbf{Slide} & : &$\Y[z_1,z_1,z_2] =  \Y[s,z_1,z_2],$  
\end{tabular}\\[0.5cm] 
where $z_0, z_1, z_2, z_3 \in G$.
For $(G,s)\rTo^f (G',s')$ a morphism  
in $\mathcal{A}b_{s}$, $\mathcal{A}_1(f)$ maps each  
generator $\Y[z_1,z_2,z_3]$ of  $\mathcal{A}_1(G,s)$ to  
$\Y[f(z_1),f(z_2),f(z_3)] \in \mathcal{A}_1(G',s')$. 
\begin{example} 
\label{ex:functors}
The map $\left[G \rMapsto (G,0)\right]$ makes 
$\mathcal{A}b$ a (full) subcategory of $\mathcal{A}b_s$.  
It follows from the definitions that the 
following diagram commutes: 
\begin{diagram} 
\mathcal{A}b & \rInto & \mathcal{A}b_s \\ 
& \rdTo<{\Lambda^3(-)} & \dTo>{\mathcal{A}_1} \\ 
& & \mathcal{A}b. 
\end{diagram} 
\end{example} 
Non-trivial examples will be given in the next paragraph. 
For future use, note that this category has an obvious  
pull-back construction extending that of $\mathcal{A}b$: 
\begin{diagram} 
(G_1,s_1)\times_{(G,s)} \SEpbk (G_2,s_2)  & \rTo & (G_2,s_2) \\ 
 \dTo &  & \dTo>{f_2} \\ 
 (G_1,s_1) & \rTo_{f_1} & (G,s) 
\end{diagram} 
where $(G_1,s_1)\times_{(G,s)} (G_2,s_2)$ is the subgroup of 
$G_1\times G_2$ consisting of those $(z_1,z_2)$ such that 
$f_1(z_1)=f_2(z_2)$, and with special element $(s_1,s_2)$. 

\subsection{Spin structures and the special Abelian group $P$} 
\label{subsec:spin}
In this paragraph, let $M$ be a compact oriented $3$-manifold
endowed with a Riemannian metric, and let
$FM$ be its bundle of oriented orthonormal frames: 
\begin{diagram}
SO(3) & \rInto^i & E(FM) & \rOnto^p & M.
\end{diagram}
Let $s \in H_1\left(E(FM);\mathbf{Z}\right)$ be the image by $i_*$
of the generator of $H_1\left(SO(3);\mathbf{Z}\right)\simeq 
\mathbf{Z}_2$.
Recall that $M$ is spinnable and that $Spin(M)$ can be defined as 
\begin{displaymath}
Spin(M)=\left\{y\in H^1\left(E(FM);\mathbf{Z}_2\right),
\ <y,s>\neq 0\right\},
\end{displaymath}
which is essentially independent of the metric. The manifold 
$M$ being spinnable, $s$ is not $0$ (and so is of order $2$).

Now, $Spin(M)$ being an affine space over $H^1(M;\mathbf{Z}_2)$ with action given by
\begin{displaymath}
\forall x \in H^1(M;\mathbf{Z}_2), \forall \sigma \in Spin(M),
\quad x\cdot \sigma:= \sigma + p^*(x),
\end{displaymath}
we can consider the space
\begin{displaymath} 
A\left( Spin(M),\mathbf{Z}_2\right) 
\end{displaymath}
of $\mathbf{Z}_2$-valued affine functions on $Spin(M)$. 
For instance, $\overline 1 \in A\left( Spin(M),\mathbf{Z}_2\right)$ 
will denote the constant map defined by $\sigma \rMapsto 1$.\\
There is a canonical map
\begin{displaymath}
A\left( Spin(M),\mathbf{Z}_2\right) \rTo^\kappa  H_1(M;\mathbf{Z}_2).
\end{displaymath}
For $f\in A\left( Spin(M),\mathbf{Z}_2\right)$, 
the homology class $\kappa(f)$ is defined unambiguously by
\begin{displaymath}
\forall \sigma, \sigma' \in Spin(M),\quad
f(\sigma')-f(\sigma)=<\sigma'/\sigma,\kappa(f)>\in \mathbf{Z}_2,
\end{displaymath}
where $\sigma'/\sigma \in H^1(M;\mathbf{Z}_2)$ is defined by the affine action
of $H^1(M;\mathbf{Z}_2)$ on $Spin(M)$. Another canonical map is
$$
H_1\left(E(FM);\mathbf{Z}\right) \rTo^{e}  A\left( Spin(M),\mathbf{Z}_2\right)
$$
sending a $x$ to the map defined by $\sigma \rMapsto <\sigma,x>$.
Next lemma gives us a nice understanding of the
special Abelian group $\left(H_1\left(E(FM);\mathbf{Z}\right),s\right)$.
\begin{lemma} 
\label{lem:type}
a) The following diagram of special groups is a pull-back diagram:
\begin{diagram}
\left(H_1\left(E(FM);\mathbf{Z}\right),\SEpbk s\right)   
& \rTo^{e} & 
\left(A\left( Spin(M),\mathbf{Z}_2\right),\overline{1}\right) \\ 
\dTo<{p_*} &  & \dTo>{\kappa} \\ 
\left(H_1\left(M;\mathbf{Z}\right),0\right) & \rTo_{-\otimes\mathbf{Z}_2} & 
\left(H_1\left(M;\mathbf{Z}_2\right),0\right). 
\end{diagram}
b) Let $t$ be the map
$$
\big\{ \textrm{Oriented framed knots in }M \big\}   
\rTo^t  H_1\left(E(FM);\mathbf{Z}\right) 
$$  
which adds to any oriented framed knot $K$ an extra $(+1)$-twist, and 
next sends it to the homology class of its lift in $FM$. Then,
\begin{enumerate}
\item[(i)] $t$ is surjective;
\item[(ii)] $t_{K_1}=t_{K_2}$ if and only if $K_1$ and $K_2$  
are cobordant as oriented knots in  
$M$ and if their framings with respect to a surface with 
boundary $(K_1)\cup(-K_2)$ then differ from each other by an even integer;
\item[(iii)] if $K_1\sharp K_2$ denotes the  band connected sum of $K_1$ 
and $K_2$, then $t_{K_1\sharp K_2}=t_{K_1}+t_{K_2}$;
\item[(iv)] the $k$-framed trivial oriented knot ($k\in \mathbf{Z}$)
is sent by $t$ to $k\cdot s$.
\end{enumerate}
\end{lemma}
\noindent \textbf{Proof.} We begin by proving a). 
The commutativity of the diagram of special groups is easy to
verify. 
By functoriality, we get a map
$$
\left(H_1\left(E(FM);\mathbf{Z}\right), s\right)  \rTo^{(p_*,e)}  
\left(H_1\left(M;\mathbf{Z}\right),0\right) 
\times_{\left(H_1\left(M;\mathbf{Z}_2\right),0\right)}
\left(A\left( Spin(M),\mathbf{Z}_2\right),\overline{1}\right).
$$
The Serre sequence associated to the fibration $FM$ gives for
homology with integer coefficients:
\begin{diagram}
0 & \rTo & H_1(SO(3);\mathbf{Z}) & \rTo^{i_*} &
H_1(E(FM);\mathbf{Z}) & \rTo^{p_*} & H_1(M;\mathbf{Z}) & \rTo & 0.
\end{diagram}
The bijectivity of $(p_*,e)$ follows from the exactness of this sequence.

We now prove b) and we begin with assertion (iv). Let $K$ be
a trivial $k$-framed oriented knot, 
let $\ast \in K$ and let $e=(e_1,e_2,e_3)\in
p^{-1}(\ast)$ be the framing of $K$ at $\ast$. We denote by $\tilde{K}$
the lift of $K$ to $FM$. Then, as a loop in $E(FM)$, $\tilde{K}$ is
homotopic to the loop in the fiber $p^{-1}(\ast)$ defined by
\begin{displaymath} 
[0,1]\ni t\rMapsto R_{2\pi(k+1)t}(e),
\end{displaymath} 
where $R_{\theta}$ (with $\theta \in \mathbf{R}$) denotes 
the rotation of oriented axis directed by $e_3$
and angle $\theta$. From an appropriate description of the
generator of $\pi_1\left(SO(3)\right)\simeq\mathbf{Z}_2$,
it follows that $\left[\tilde{K}\right]=(k+1)\cdot s \in 
H_1(E(FM);\mathbf{Z})$, and assertion (iv) then follows.

Let us make an observation. Let $K$ be any oriented framed knot in $M$;
since the framing  of $K$ determines a 
trivialization of its normal bundle in $M$, it allows
us to restrict any spin structure on $M$ to $K$. Recall now that 
the cobordism group  $\Omega_1^{Spin}$ is isomorphic to $\mathbf{Z}_2$
(with generator given by $\mathbf{S}^1$ endowed with the spin structure
induced by its Lie group structure:
see \cite[p. 35, 36]{Kirby}).
The following observation then makes sense: 
\begin{equation}
\label{eq:Boolean}
\forall \sigma\in Spin(M), \quad
e(t_K)(\sigma)= \left(K,\sigma|_K\right) \ \in
\Omega_1^{Spin}\simeq \mathbf{Z}_2,
\end{equation}  
and can be derived from an appropriate characterization of the spin structures 
on the circle (see \cite[p. 35, 36]{Kirby}).

Let now $K_1$ and $K_2$ be some disjoint oriented framed knots 
in $M$. There is an obvious genus $0$ surface with boundary 
$K_1\sharp K_2 \dot\cup (-K_1)\dot\cup (-K_2)$. Then, according
to (\ref{eq:Boolean}), we have $e(t_{K_1\sharp K_2})=
e(t_{K_1})+e(t_{K_2})$. Also, $p_*\left( t_{K_1\sharp K_2}\right)
$$=[K_1\sharp K_2]$$=[K_1]+[K_2]$$=p_*\left(t_{K_1}\right)+
p_*\left(t_{K_2}\right)$, and so by a), we obtain
that assertion (iii) holds for $K_1$ and $K_2$.\\
We now justify assertion (ii).
According to a), $t_{K_1}=t_{K_2}$ if and only if
$p_*(t_{K_1})=p_*(t_{K_2})$ and $e(t_{K_1})=e(t_{K_2})$.
Also, the condition $p_*(t_{K_1})=p_*(t_{K_2})$ holds
if and only if $K_1$ and $K_2$ are homologous in $M$.
In this case, let $S$ be an embedded oriented surface in $M$ 
such that $\partial S=K_1\dot \cup (-K_2)$.
Let $k_i$ be the framing of $K_i$ with respect to $S$ and let $K'_i$
be the oriented framed knot obtained from $K_i$ by adding
an extra $(-k_i)$-twist, so that the framing of $K'_i$ is
given by $S$. Then, according 
to (\ref{eq:Boolean}), we have $e\left(t_{K'_1}\right)=e\left(t_{K'_2}\right)$.
Moreover, applying assertions (iii) and (iv), we obtain:
$e\left(t_{K'_i}\right)=e\left(t_{K_i}\right)+k_i\cdot s$.
We conclude that $e\left(t_{K_1}\right)=e\left(t_{K_2}\right)$
if and only if $k_1$ and $k_2$ are equal modulo $2$, proving thus
assertion (ii).
 
Let $x \in H_1(E(FM);\mathbf{Z})$,
then $p_*(x) \in H_1(M;\mathbf{Z})$ can be realized by an oriented 
knot $K$ in $M$: we give it an arbitrary framing. By construction,
$p_*(t_K - x)=0 \in H_1(M;\mathbf{Z})$, and so by exactness of the
Serre sequence, $t_K - x=\varepsilon \cdot s$ with $\varepsilon\in \{0,1\}$.
By possibly band-summing $K$ with a trivial
$(+1)$-framed knot when $\varepsilon=1$, 
and according to assertion (iii) and (iv), 
the framed knot $K$ can be supposed to be 
such that $t_{K}=x$; this proves assertion (i).
$\square$\\

We now restrict ourselves to the $3$-manifold $M=1_\Sigma=\Sigma\times I$
where $\Sigma$ can be $\Sigma_g$ or $\Sigma_{g,1}$.
The inclusion $i^+:\Sigma \rInclus 1_{\Sigma}$,
with image $\Sigma^+$, induces an isomorphism
between $H$ and $H_1(M;\mathbf{Z})$ and a bijection between
$Spin(\Sigma)$ and $Spin(M)$. As shown by Johnson in \cite{JSpin},
there is an algebraic way to think of $Spin(\Sigma)$. Indeed, there exists
a canonical affine isomorphism
$$
Spin(\Sigma)  \rTo^{\simeq}  \Omega_g,
$$
sending any spin structure $\sigma$ to a quadratic form $q_{\sigma}$ 
which can be defined as follows.
Let $x\in H_{(2)}=H_1(\Sigma;\mathbf{Z}_2)$ be represented by an oriented simple
closed curve on $\Sigma^+$; by framing it along $\Sigma^+$ and pushing it into the interior of
$1_\Sigma$, we get a framed oriented knot $K$ in $1_{\Sigma}$. Then,
\begin{equation}
\label{eq:def_quad}
q_\sigma(x)=e(t_K)\left(\sigma \times I\right)\ \in \mathbf{Z}_2.
\end{equation}
Therefore, according to Lemma \ref{lem:type} a), 
$\left(H_1\left(E\left(F1_\Sigma\right);\mathbf{Z}\right),
s\right)$ is canonically isomorphic to the special Abelian group
defined by the pull-back construction
\begin{diagram} 
(H,0) \times_{\left(H_{(2)},0\right)}  \SEpbk 
\left(B^{(1)}_g,\overline{1}\right)
  & \rTo^e & \left(B^{(1)}_g,\overline{1}\right) \\ 
\dTo<p &  & \dTo>{\kappa} \\ 
(H,0)  & \rTo_{-\otimes\mathbf{Z}_2} & \left(H_{(2)},0\right)  
\end{diagram}
whose projections are denoted by $p$ and $e$, and where
$\kappa$ is the composite
\begin{diagram}
B^{(1)}_g & \rOnto & B^{(1)}_g/B^{(0)}_g &
\rTo^\simeq & H_{(2)}.
\end{diagram}
The last isomorphism here identifies $\overline{h}$ with 
$h_{(2)}$ for all $h\in H$ (this is well-defined by 
(\ref{eq:quadratic_identity}) and (\ref{eq:iso_Boolean})).
We define the special Abelian group $P$ to be
\begin{displaymath}
\begin{array}{|c|}
\hline \\
P=(H,0) \times_{\left(H_{(2)},0\right)}  
\left(B^{(1)}_g,\overline{1}\right), \\
\\
\hline
\end{array}
\end{displaymath}
and $\mathcal{A}_1(P)$ is the space of graphs announced in the introduction.
\begin{remark} \label{rem:elts_of_P}
Thus, any element $z$ of $P$ can be written as
\begin{displaymath}
z=\left(h,\overline{h}+\varepsilon\cdot \overline{1}\right) \in P, 
\end{displaymath}
with $h \in H$ and $\varepsilon \in \{0,1\}$. Observe also the following.
Suppose that there exists a simple oriented closed curve in $\Sigma^+$ with homology class $h$. 
Let $K$ be the push-off of this curve, framed along $\Sigma^+$, 
with an extra $\varepsilon$-twist. Then, it follows from (\ref{eq:def_quad}) that
$t_K=z \in P \simeq H_1\left(E\left(F1_\Sigma\right);\mathbf{Z}\right)$.
\end{remark} 
\begin{remark} \label{rem:ext}
According to the proof of Lemma \ref{lem:type},
the Serre sequence for homology associated to the bundle $F1_\Sigma$
gives the following short exact sequence:
\begin{diagram} 
0 & \rTo & \mathbf{Z}_2 & \rTo & P & \rTo^p & H & \rTo 0, 
\end{diagram}
where  $\mathbf{Z}_2$ is injected into $P$ by sending $1$
to $(0,\overline 1)$. The map $s:H\rTo P$ defined by 
$s(h)=\left(h,\overline{h}\right)$ is a set-theoretic section. According to 
(\ref{eq:quadratic_identity}), the associated $2$-cocycle
$H\times H\rTo \mathbf{Z}_2$ is the mod $2$ reduced
intersection form of $\Sigma$. Thus, $P$ 
is isomorphic to $H\rtimes \mathbf{Z}_2$ with crossed product defined by 
\begin{displaymath} 
(h_{1},\varepsilon_{1}) \cdot (h_{2},\varepsilon_{2}) =  
      (h_{1}+h_{2},\varepsilon_{1}+ \varepsilon_{2}+ h_{1}\bullet h_{2}). 
\end{displaymath} 
The element $\left(h,\overline{h}+\varepsilon\cdot\overline{1}\right)\in P$ 
corresponds to $(h,\epsilon)\in H\rtimes \mathbf{Z}_2$.  
\end{remark}

\subsection{The surgery map $\psi_1$}  
\begin{figure}[!h] 
\begin{center} 
\includegraphics[width=7.5cm,height=6cm]{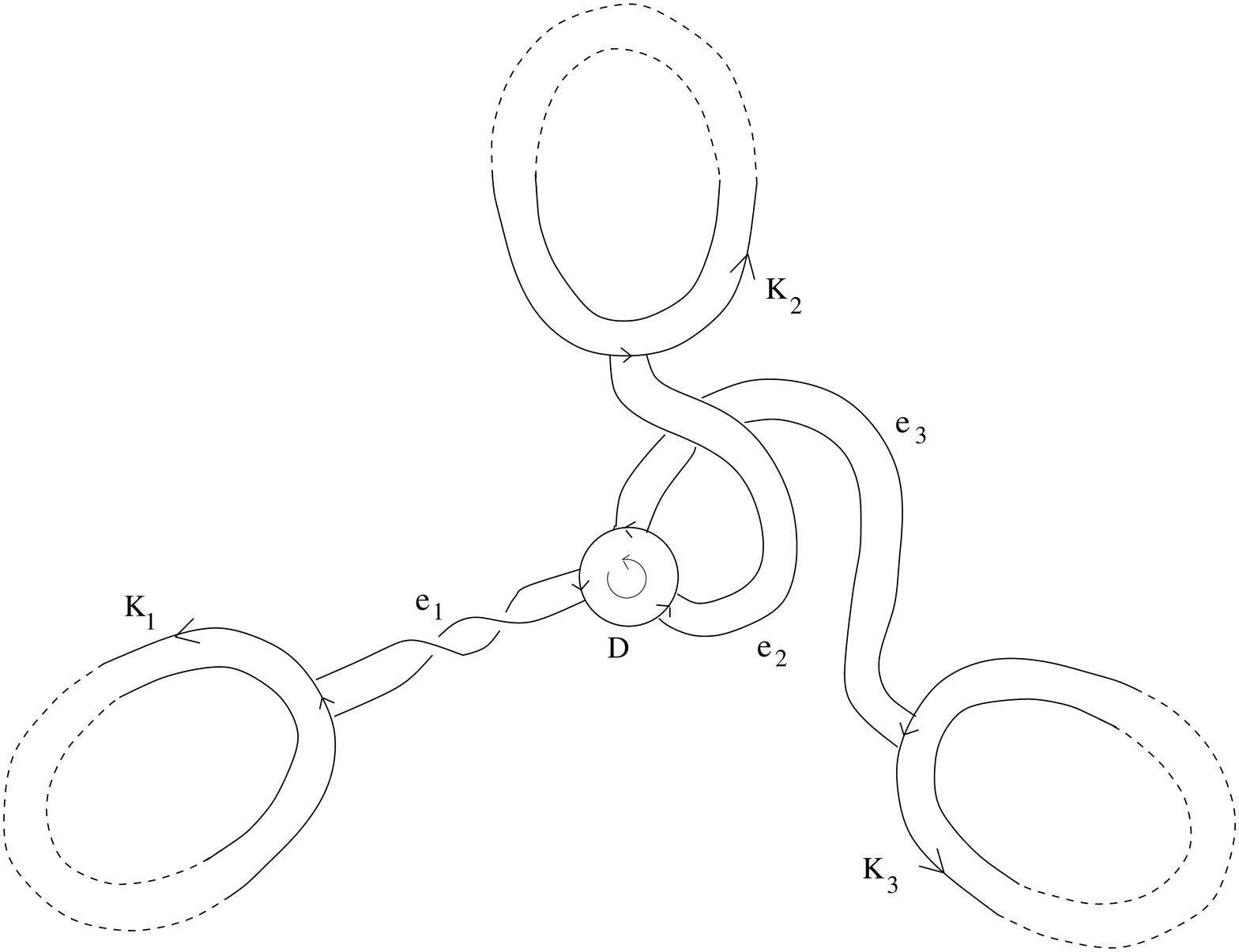} 
\fcaption{Embedding the Y-graph} 
\label{fig:embedding} 
\end{center} 
\end{figure}

In this paragraph,  $\Sigma$ is allowed to be $\Sigma_g$ or 
$\Sigma_{g,1}$ and the surgery map 
$\psi_1:\mathcal{A}_1(P) \rTo \overline{\mathcal{C}}_1(\Sigma)$ 
is constructed by means of calculi of clovers.
\label{subsec:psi1} 
\begin{conventions}
\label{conv:Y}
Here, we \emph{adopt} Goussarov's convention for the surgery meaning 
of $Y$-graphs and clovers \cite{G}, \cite{GGP}.  
\end{conventions}
\noindent
Denote by $\tilde{\mathcal A}_1(P)$ the free Abelian group generated  
by abstract Y-shaped graphs whose univalent vertices are labelled  
by $P$, and which are equipped with an orientation at their trivalent vertex:  
$\mathcal{A}_1(P)$ is a quotient of $\tilde{\mathcal A}_1(P)$. 
For each generator $\mathsf{Y}[z_1,z_2,z_3]$ of $\tilde{\mathcal A}_1(P)$, 
where $z_i\in P$, pick 
some disjoint oriented framed knots $K_i$  
in the interior of $1_\Sigma$ such that  
$t_{K_i}=z_i \in P\simeq 
H_1\left(E\left(F 1_\Sigma\right);\mathbf{Z}\right)$; 
this is possible according to Lemma \ref{lem:type} b) (i). 
Next, pick an embedded $2$-disk $D$ in the interior 
of $1_\Sigma$ and disjoint from the $K_i$'s,  
orient it in an arbitrary way, and connect it to the $K_i$'s with some 
bands $e_i$. These band sums are required to be compatible  
with the orientations, and to be coherent  
with the cyclic ordering $(1,2,3)$. See Fig. \ref{fig:embedding}
as an illustration. 
What we obtain in $1_\Sigma$ is precisely a \emph{Y-graph},  
as defined by Goussarov in \cite{G}. We denote it by  
$\phi\big(\mathsf{Y}[z_1,z_2,z_3]\big)$. For example, 
as follows from Lemma \ref{lem:type} b) (iv), if $z_1$ 
is the special element $s$ of $P$, the corresponding 
leaf $K_1$ of $\phi\big(\mathsf{Y}[s,z_2,z_3]\big)$  
can be chosen to be unknotted and $(+1)$-framed; such a leaf is  
called \emph{special} in \cite{GGP}.\\ 
We now put  
$\tilde{\psi}_1\big(\mathsf{Y}[z_1,z_2,z_3]\big)$ 
to be the $Y_2$-equivalence class of the surgered 
manifold $(1_{\Sigma})_{\phi(\mathsf{Y}[z_1,z_2,z_3])}$, 
so that  we get a map  
$$
\tilde{\mathcal{A}}_1(P) \rTo^{\tilde{\psi}_1} \overline{\mathcal{C}}_1(\Sigma). 
$$
\begin{theorem} 
\label{th:surgery_map} 
The map $\tilde{\psi}_1$ does not depend on the choice of $\phi$, and induces 
a surjective quotient map 
$$
\mathcal{A}_1(P)  \rOnto^{\psi_1}  \overline{\mathcal{C}}_1(\Sigma). 
$$
\end{theorem} 
\noindent \textbf{Proof.} The proof might be read with a copy of \cite{GGP} in hand. 
Using the above notation, we begin with showing that 
$\tilde{\psi}_1\big(\mathsf{Y}[z_1,z_2,z_3]\big)$ does not depend on the choice  
of $\phi\big(\mathsf{Y}[z_1,z_2,z_3]\big)$. For this, we recall two facts 
concerning any Y-graph $G$ in a homology cylinder $M$ (see Remark \ref{rem:facts} below): 
\begin{description} 
\item[Fact 1] the $Y_2$-equivalence class of $M_G$ is not modified  
when an edge of $G$ is band-summed  
with a (disjoint) oriented framed knot of $M$; 
\item[Fact 2] the $Y_2$-equivalence class of $M_G$  
is inverted when an edge of $G$ is half-twisted. 
\end{description} 
Using these, the independance on the choice of the disk $D$, its orientation 
and the edges $e_i$ is easily shown.\\ 
We now show the independance on the choice of the leaves $K_i$. 
Suppose for example that $K'_1$ is another choice of $K_1$. Then, 
according to Lemma \ref{lem:type} b) (ii), there  exists an embedded  
oriented surface $F$ in $1_\Sigma$ such that  
$\partial F=K_1\dot \cup (-K'_1)$ and such that, if $k$ (resp. $k'$) 
is the framing of $K_1$ (resp. $K'_1$) with respect to $F$, 
$(k-k')$ is even. We also assume transversality of $F$ with the edges 
of the Y-graph, and with the two other leaves $K_2$ and $K_3$. 
Let $g(F)$ denote the genus of $F$, let 
$m$ be the number of intersection points 
of $F$ with the edges, and for $i=2,3$, let $n_i$ be the number of 
intersection points of $F$ with $K_i$. 
If all of the integers $g(F)$, $(k-k')$, $m$, $n_2$ and $n_3$ are zero,    
the two Y-graphs are isotopic and we are done. 
In the general case, recall from \cite[\S4.3]{GGP} that there is 
a procedure for \emph{simplifying the leaves}. 
The main tool for this is the following:   
\begin{description}   
\item[Fact 3] if $G_1$ and $G_2$ are two Y-graphs in $1_\Sigma$ 
obtained from a Y-graph $G$ by \emph{splitting a leaf}, then  
$(1_{\Sigma})_G=(1_{\Sigma})_{G_1}\cdot (1_{\Sigma})_{G_2}  
\in \overline{\mathcal{C}}_{1}(\Sigma)$ (see Remark \ref{rem:facts}).  
\end{description} 
Splitting $\big(g(F)+|k-k'|/2+m+n_2+n_3\big)$ times the leaf $K_1$,  
splitting  $n_2$ times the leaf $K_2$ and splitting $n_3$ times  
the leaf $K_3$, we see that 
the result $\tilde{\psi}_1\big(\mathsf{Y}[z_1,z_2,z_3]\big)$  
in $\overline{\mathcal{C}}_1(\Sigma)$  defined  
by the choice of $K_1$ differs 
from the one defined by $K'_1$ 
by some elements of the form $(1_\Sigma)_G$, 
where $G$ satisfies one of the following conditions: 
\begin{enumerate} 
\item[(i)] $G$ has a leaf which bounds a genus $1$ surface disjoint from $G$ 
and with respect to which the leaf is $0$-framed; 
\item[(ii)] $G$ has a leaf which bounds a disk disjoint from $G$, and with respect 
to which the leaf is $(\pm 2)$-framed; 
\item[(iii)] $G$ has a leaf which bounds a disk with respect to which it 
is $0$-framed, and this disk intersects $G$ in exactly one point belonging 
to an edge;   
\item[(iv)] $G$ has two leaves which are linked as the Hopf link. 
\end{enumerate} 
Let us now verify that all of these elements  
vanish in $\overline{\mathcal{C}}_1(\Sigma)$. If $G$ is of type (i), 
the surgery effect of $G$ is the same as a clover of degree $2$ 
(apply \cite[Lem. 5.1]{GGP} and \cite[Th. 2.4]{GGP}).  
If $G$ is of type (ii), by again cutting its leaf we get 
$(1_\Sigma)_G=2\cdot(1_\Sigma)_{G'}$ where $G'$ has a special leaf; but  
$(1_\Sigma)_{G'}=-(1_\Sigma)_{G'}$ by Fact 2. 
If $G$ is of type (iii), by applying Fact 1 the edge can be slid away    
from the leaf, we then get a Y-graph with a \emph{trivial} leaf  
which has no surgery effect by the ``blow-up move'' of \cite[Fig. 6]{GGP}.
If $G$ is of type (iv), 
by applying \cite[Th. 2.4]{GGP}, we obtain a Y-graph  
with a looped edge, but this is stated to be $0$ in  
$\overline{\mathcal{C}}_1(\Sigma)$ by the so-called 
\emph{LOOP relation}. This relation is easily shown  
from \cite[Lem. 2.3]{GGP} and from Fact 1 and Fact 2. 
This completes the proof of the independance of $\tilde{\psi}_1$ on $\phi$.
 
The fact that $\tilde{\psi}_1$ is surjective follows immediately from the fact that the Abelian group  
$\overline{\mathcal{C}}_{1}(\Sigma)$ is generated by the  
homology cylinders $(1_\Sigma)_{G}$ where $G$ is a single 
Y-graph (this is also proved by standard calculi of clovers). 
 
We now show that the map $\tilde{\psi}_1$ factors 
through $\tilde{\mathcal{A}}_1(P) \rTo \mathcal{A}_1(P)$. 
The AS relation is proved in $\overline{\mathcal{C}}_{1}(\Sigma)$ 
from Fact 2 and an isotopy of the $Y$-graph -- see \cite[Cor. 4.6]{GGP}.\\ 
The multilinearity relation follows from Fact 3. 
Indeed, let $G$ be a Y-graph in $1_\Sigma$ with $K$ as a leaf. Split 
the leaf $K$ to $K_1$ and $K_2$, and let $G_1$ and $G_2$ be the corresponding 
new Y-graphs. Then, $(1_\Sigma)_G=(1_\Sigma)_{G_1}\cdot (1_\Sigma)_{G_2}  
\in \overline{\mathcal C}_1(\Sigma)$.  
Since $K$ is the band connected sum of $K_1$ and $K_2$, we have by Lemma 
\ref{lem:type} b) (iii): $t_K=t_{K_1}+t_{K_2} \in P$.\\ 
The slide relation is shown to be satisfied in  
$\overline{\mathcal{C}}_{1}(\Sigma)$ thanks to the ``leaf slide'' 
move of \cite[Fig. 6]{GGP}. For this, let $G$ be a Y-graph  
in $1_\Sigma$ with some leaves $K_1$ and $K_2$ such that  
$t_{K_1}=-t_{K_2}$. By sliding the leaf $K_2$ along $K_1$, 
we obtain a new Y-graph $G'$ with the same surgery effect as $G$, 
such that $K'_1=K_1$ and such that $K'_2$ is the band 
connected sum of $K_1$ 
and $K_2$ with an extra $(-1)$-twist. So, by Lemma  
\ref{lem:type} b) (iii) and (iv), we have $t_{K'_2}=t_{K_1}+t_{K_2}+s=s\in P.$
This shows that the relation $\mathsf{Y}[z_1,-z_1,z_3] =  
\mathsf{Y}[z_1,s,z_3]$ ($z_1,z_3 \in P$)  
is satisfied in $\overline{\mathcal{C}}_{1}(\Sigma)$. 
The slide relation, as stated in \S 2.1, follows 
then from the AS and multilinearity relations. 
$\square$\\

\begin{remark}
\label{rem:facts}
The proof of Fact 1, Fact 2 and Fact 3 use calculus of clovers and can respectively 
be obtained from the proof of Cor. 4.2, Lem. 4.4 and Cor. 4.3 in \cite{GGP}. 
Alternatively, those facts can be considered as corollaries of these results in the following way. 
Denote by $\mathbf{Z} \mathcal{C}_1\left(\Sigma\right)$ the free Abelian group generated by
the set $\mathcal{C}_1\left(\Sigma\right)$, and let 
$$
\mathbf{Z} \mathcal{C}_1\left(\Sigma\right)=\mathcal{F}^Y_0\left(1_\Sigma\right) \supset
\mathcal{F}^Y_1\left(1_\Sigma\right)\supset \mathcal{F}^Y_2\left(1_\Sigma\right) \supset \cdots 
$$
be its Goussarov-Habiro filtration \cite[\S 1.4]{GGP}. Results are 
stated in \cite{GGP} to hold in the graded space $\mathcal{G}_k\left(1_{\Sigma}\right)=
\mathcal{F}^Y_k\left(1_\Sigma\right)/\mathcal{F}^Y_{k+1}\left(1_\Sigma\right)$. 
Consider also the homomorphism of Abelian groups 
$$
\mathbf{Z} \mathcal{C}_1\left(\Sigma\right)\rTo^\upsilon
\overline{\mathcal{C}}_1\left(\Sigma\right)
$$
which assigns to any homology cylinder its $Y_2$-equivalence class. 
The invariant $\upsilon$ is primitive, in the sense that it restricts to 
$\mathcal{C}_1\left(\Sigma\right)$ to a monoid homomorphism, and is a degree $1$
invariant\footnote{In fact, $\upsilon$ is a universal degree 1 primitive invariant 
for homology cylinders. See \cite[\S 6.4]{H} for a similar invariant, of any degree, for knots
in the $3$-sphere.}\,\, as follows from calculus of clovers.
In particular, $\upsilon$ induces a homomorphism $\mathcal{G}_1\left(1_{\Sigma}\right)\rTo
\overline{\mathcal{C}}_1\left(\Sigma\right)$, by which Fact 1, Fact 2 and Fact 3 are
respectively the images of Cor. 4.2, Lem. 4.4 and Cor. 4.3.
\end{remark}
\section{Johnson Homomorphism and  Birman-Craggs Homomorphisms for  
Homology Cylinders} 
\label{sec:homomorphisms} 
In this section, the first Johnson homomorphism and the Birman-Craggs 
homomorphisms are extended to the monoid of homology cylinders.

\subsection{The first Johnson homomorphism for homology cylinders} 
In \cite{GL} the notion of Johnson homomorphisms for 
homology cobordisms over $\Sigma_{g,1}$ was introduced.
In this paragraph,
we allow $\Sigma$ to be $\Sigma_g$ or $\Sigma_{g,1}$, and give
the definition of the first Johnson homomorphism in both cases.

\begin{figure}[!h] 
\begin{center} 
\includegraphics[width=7cm,height=3.5cm]{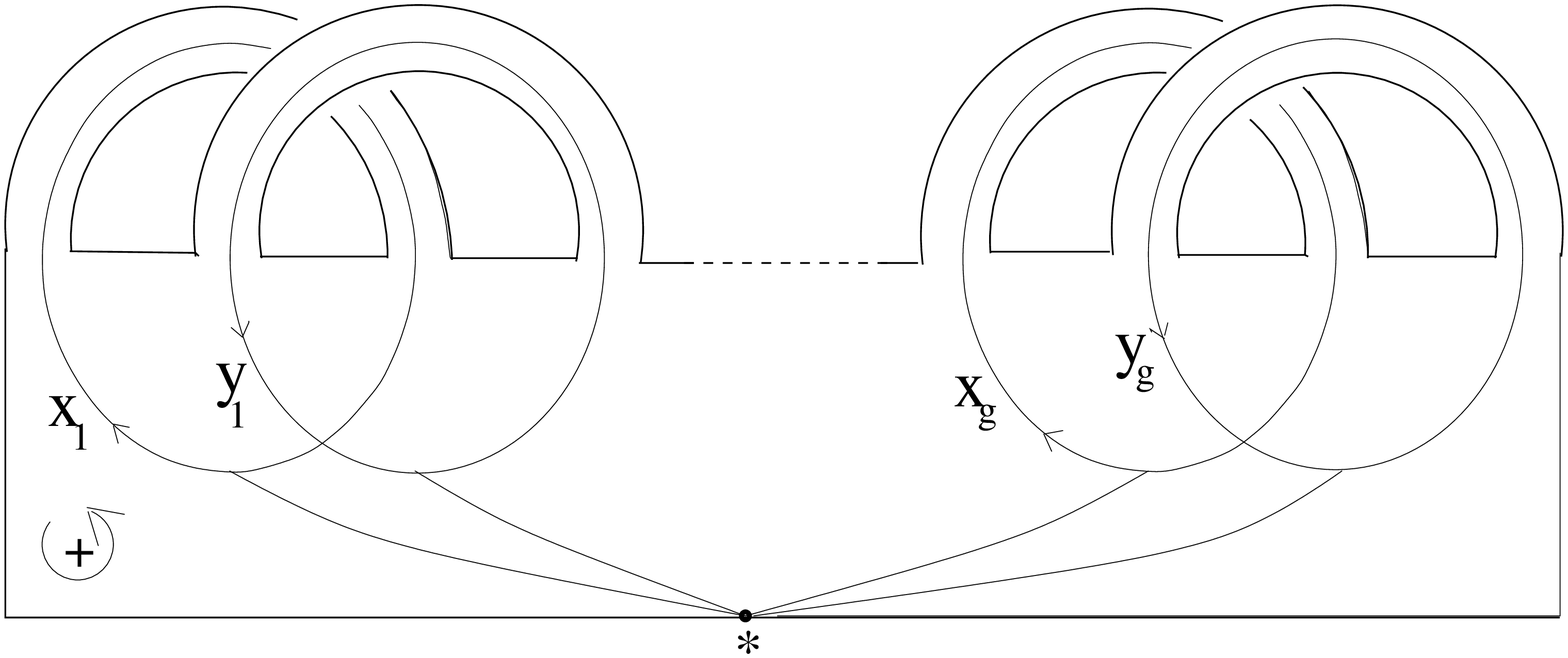}
\fcaption{The based curves $(x_i,y_i)_{i=1}^g$ on $\Sigma_{g,1}$} 
\label{fig:basis} 
\end{center} 
\end{figure}

The fundamental group of $\Sigma$ with base point $\ast\in \Sigma$ will 
be denoted by $\pi^{(\ast)}$, 
and $\pi^{(\ast)}_{k}$ will denote the 
$k^{th}$ term of its lower central 
series, beginning at $\pi^{(\ast)}_1=\pi^{(\ast)}$.
We denote by $(x_i,y_i)_{i=1}^g$ the based loops
depicted in Fig. \ref{fig:basis}
or their corresponding images under an inclusion 
$\Sigma_{g,1}\subset \Sigma_g$. Then,
\begin{displaymath} 
\begin{array}{ll}
\textrm{in the boundary case, } & 
\pi^{(\ast)}=F(x_1,\dots,x_g,y_1,\dots,y_g) , \\
&\\
\textrm{and in the closed case, } & \pi^{(\ast)}=
\langle x_1,\dots,x_g,y_1,\dots,y_g \ | \prod_{i=1}^g[x_i,y_i]=1 \rangle.
\end{array}
\end{displaymath} 

Given a homology cobordism $(M,i^{+},i^{-}) 
\in \mathcal{C}(\Sigma)$, the map
$i^{\pm}$ induces an isomorphism at the level of each nilpotent 
quotient (by Stallings \cite{S}). 
We choose a path $\gamma \subset M$ going from
$i^+(\ast)$ to $i^-(\ast)$, and then consider the following composite:
\begin{diagram}
\frac{\pi^{(\ast)}}{\pi^{(\ast)}_3} & \rTo^{i^+_3}_\simeq &
\frac{\pi_1(M,i^+(\ast))}{\pi_1(M,i^+(\ast))_3} & \rTo^{c_\gamma}_\simeq &
\frac{\pi_1(M,i^-(\ast))}{\pi_1(M,i^-(\ast))_3} &
\rTo^{(i^-_3)^{-1}}_\simeq & \frac{\pi^{(\ast)}}{\pi^{(\ast)}_3}.
\end{diagram}
Up to inner automorphisms, this is independent on the choice of $\gamma$,
so that there is a well-defined map
\begin{diagram}
\mathcal{C}(\Sigma) & \rTo^{\eta^{(\ast)}_1} & 
Out\left(\frac{\pi^{(\ast)}}{\pi^{(\ast)}_3}\right), 
\end{diagram}
satisfying $\eta^{(\ast)}_1(M\cdot N)=\eta^{(\ast)}_1(N)\cdot
\eta^{(\ast)}_1(M)$.
Let $\star$ be another base point in $\Sigma$, and $\gamma$ an 
arbitrary path between $\ast$ and $\star$. Conjugation by $\gamma$ induces an 
isomorphism $Out\left(\pi^{(\ast)}/\pi^{(\ast)}_3\right)\simeq 
Out\left(\pi^{(\star)}/\pi^{(\star)}_3\right)$. This 
isomorphism is independent on the choice of the path $\gamma$, and 
the maps $\eta^{(\ast)}_1$ and $\eta^{(\star)}_1$ are compatible
through it. Therefore, we get a well-defined group denoted by
$Out(\pi/\pi_3)$ and an anti-homomorphism of monoids
\begin{equation}
\label{eq:brut}
\mathcal{C}(\Sigma) \rTo^{\eta_1} 
Out\left(\frac{\pi}{\pi_3}\right). 
\end{equation} 
If we restrict ourselves to homology cylinders, we are led to a map
\begin{diagram}
\mathcal{C}_1(\Sigma) & \rTo^{\eta_1} &  Ker\left(Out\left(\frac{\pi}
{\pi_3}\right) \rightarrow Out\left(\frac{\pi}{\pi_2}\right)\right).
\end{diagram}
Observe the following exact sequence:
\begin{diagram} 
1 & \rTo & Hom\left(H,\frac{\pi^{(\ast)}_2}{\pi^{(\ast)}_3}\right) & \rTo &  
Aut\left(\frac{\pi^{(\ast)}}{\pi^{(\ast)}_3}\right)  
 & \rTo & Aut\left(\frac{\pi^{(\ast)}}{\pi^{(\ast)}_2}\right)
\end{diagram}
where any $f\in 
Hom\left(H,\pi^{(\ast)}_2/\pi^{(\ast)}_3\right)$ is sent
to the automorphism of $\pi^{(\ast)}/\pi^{(\ast)}_3$ 
which sends $\overline{x}$ to $\overline{x} f(\overline{x})$ (with
$x\in \pi^{(\ast)}$). Hence we have the following exact sequence:
\begin{diagram}
1 & \rTo & \frac{Hom\left(H,\pi^{(\ast)}_2/\pi^{(\ast)}_3\right)}
{[H,-]} & \rTo &  
Out\left(\frac{\pi}{\pi_3}\right) & \rTo & Out\left(\frac{\pi}{\pi_2}\right).
\end{diagram}
Here, $[H,-]$ stands for the subgroup of $Hom\left(H,\pi^{(\ast)}_2/\pi^{(\ast)}_3\right)$
consisting of those homomorphisms $[h,-]$ defined for any $h\in H$ by $x\rMapsto$$[h,x]$,
where $H$ is identified with $\pi_1^{(\ast)} /\pi_2^{(\ast)}$.
Consequently, we have defined an anti-homomorphism of monoids
\begin{displaymath}
\mathcal{C}_1(\Sigma) \rTo^{\eta_1}  
\frac{Hom\left(H,\pi^{(\ast)}_2/
\pi^{(\ast)}_3\right)}{[H,-]}.
\end{displaymath}  
In the sequel, we denote by $\mathsf{L}(H)=\oplus_n \mathsf{L}_n(H)$, 
the free Lie 
$\mathbf{Z}$-algebra 
on the $\mathbf{Z}$-module $H$, and distinguish
the boundary case from the closed case.\\

In the boundary case, as $\pi^{(\ast)}$ is free and $H$ is the Abelianized of $\pi^{(\ast)}$,
$\mathsf{L}_2(H)$ is canonically isomorphic to $\pi^{(\ast)}_2/\pi^{(\ast)}_3$.
Also, there is a sequence of isomorphisms
$Hom\left(H, \mathsf{L}_2(H)\right)\simeq H^*\otimes \mathsf{L}_2(H)
\simeq H\otimes \mathsf{L}_2(H)$, with last one induced 
by $\bullet$-duality. Through these, $[H,-]\subset
Hom\left(H, \mathsf{L}_2(H)\right)$ becomes 
$A_{g,1}\subset  H\otimes \mathsf{L}_2(H)$ defined by
\begin{displaymath}
A_{g,1}=\left\{\sum^{g}_{i=1}\left(x_i\otimes [h,y_i]-
y_i\otimes [h,x_i]\right) \ \arrowvert \ h\in H \right\}. 
\end{displaymath}
Thus, $\eta_1$ takes values in
\begin{displaymath}
\frac{Hom\left(H,\pi^{(\ast)}_2
/\pi^{(\ast)}_3\right)}{[H,-]} \simeq
\frac{H\otimes \mathsf{L}_2(H)}{A_{g,1}}.
\end{displaymath}
The group $\Lambda^3 H$ can be seen 
as a subgroup of $H\otimes \mathsf{L}_2(H)$ in the following manner:
\begin{diagram} 
0 & \rTo & \Lambda^3H & \rTo^{\nu} & H\otimes \mathsf{L}_2(H) 
& \rTo^{[-,-]} & \mathsf{L}_3(H), 
\end{diagram}
where $\nu$ is defined by 
$\nu(x\we y\we z) = x\otimes [y,z] + y\otimes [z,x]+ z\otimes [x,y]$.
Composing $\nu$ with the projection 
$H\otimes \mathsf{L}_2(H) \rOnto H\otimes \mathsf{L}_2(H)/A_{g,1}$
still gives an injection
\begin{diagram} 
\Lambda^3H & \rInto^\nu & \frac{H\otimes \mathsf{L}_2(H)}{A_{g,1}}. 
\end{diagram}
This follows from the fact that
\begin{equation} \label{eq:calcul}
\forall h\in H, \quad [h,\omega]=0 \in
\mathsf{L}_3(H) \ \Longrightarrow \ h=0,
\end{equation}
where $\omega=\sum_i [x_i,y_i] \in \mathsf{L}_2(H)$ corresponds 
\emph{via} the canonical isomorphism $\mathsf{L}_2(H)\simeq \Lambda^2H$
to the symplectic element $\omega$, defined in the introduction.\\
We now prove that $\eta_1$ takes values 
in the subgroup $\Lambda^3H$.
Suppose for this that $f\in Hom\left(H,\pi^{(\ast)}_2/
\pi^{(\ast)}_3\right) \subset
 Aut\left( \pi^{(\ast)}/\pi^{(\ast)}_3\right)$ is such that
there exists a lift $\tilde{f}\in End(\pi^{(\ast)})$ 
of $f$ fixing the boundary
element $\partial:=\prod_{i=1}^g[x_i,y_i]$ modulo
$\pi^{(\ast)}_4$. Note that this property
is verified by a representative for $\eta_1(M)$ if $M$ is a 
homology cylinder, so that proving that $f \in Ker([-,-])$ will prove
that $Im(\eta_1)\subset \Lambda^3 H$.
Let $X_i=x_i^{-1}\tilde{f}(x_i)\in \pi^{(\ast)}_2$
and $Y_i=y_i^{-1}\tilde{f}(y_i)\in \pi^{(\ast)}_2$. We have
\begin{displaymath}
\begin{array}{rcll}
   \tilde{f}(\partial) & = & \prod_i [\tilde{f}(x_i),\tilde{f}(y_i)] &\\ 
    & \equiv & \prod_i [x_i X_i,y_iY_i] &\\ 
    & \equiv & \prod_i [x_i,y_i][X_i,y_i][x_i,
      Y_i] & \textrm{ mod }\pi^{(\ast)}_4, 
\end{array}
\end{displaymath}
which implies that $\prod_i[X_i,y_i][x_i,Y_i]\equiv 1$ mod
$\pi^{(\ast)}_4$. Consequently,
\begin{displaymath}
\sum_i \left(x_i\otimes Y_i - y_i\otimes X_i\right)
\in H\otimes  \mathsf{L}_2(H),
\end{displaymath}
which essentially corresponds to $f$, goes to $0$ by the bracketting map.\\

Let us now focus on the closed case.
The canonical map $\mathsf{L}_2(H)\rOnto \pi^{(\ast)}_2/\pi^{(\ast)}_3$ 
induces an isomorphism between $\pi^{(\ast)}_2/\pi^{(\ast)}_3$
and $\mathsf{L}_2(H)/\omega$. Thus, in this case, $\eta_1$ takes values in
\begin{displaymath}
\frac{Hom\left(H,\pi^{(\ast)}_2
/\pi^{(\ast)}_3\right)}{[H,-]} \simeq
\frac{H\otimes \mathsf{L}_2(H)}{A_g}
\end{displaymath}
where $A_g=A_{g,1} + H\otimes \omega$.  
Since $\nu\left(\omega\wedge H\right)\subset A_g$, $\nu$ factors to give
\begin{diagram}
\frac{\Lambda^3H}{\omega\wedge H} & \rTo^\nu & 
\frac{H\otimes \mathsf{L}_2(H)}{A_g}.
\end{diagram}
It also follows from (\ref{eq:calcul}) that this new $\nu$ is still 
injective. Then, $\Lambda^3H/\omega\wedge H$
can be seen as a subgroup of $Hom\left(H,\pi^{(\ast)}_2/
\pi^{(\ast)}_3\right)/[H,-]$. 
Similarly to the boundary case, one shows that $\eta_1$ takes
values in $\Lambda^3H/ \omega\wedge H$.\\ 

\begin{figure}[!h] 
\begin{center} 
\includegraphics[width=12cm,height=4cm]{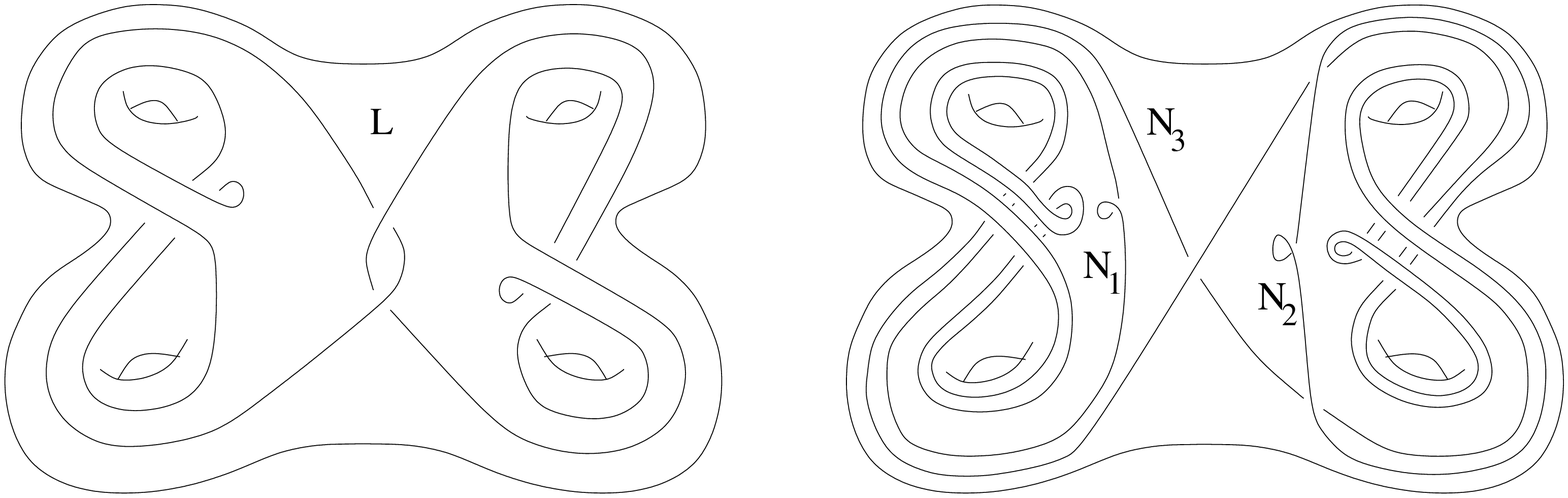}
\vspace{0.3cm} 
\fcaption{The $2$-component framed link $L$ 
and a Kirby-equivalent boundary link $N$} 
\label{fig:link} 
\end{center} 
\end{figure}

So far, we have defined some anti-homomorphisms of monoids
\begin{diagram}
\mathcal{C}_1(\Sigma_{g,1})& \rTo^{\eta_1} & \Lambda^3H &
\textrm{ and } & 
\mathcal{C}_1(\Sigma_{g})& \rTo^{\eta_1} & \frac{\Lambda^3H}{\omega\wedge H},
\end{diagram}
but next lemma allows us to go a bit further.
\begin{lemma} \label{lem:Y2} 
Let $(M,K)$ be a homology cylinder over $\Sigma$ together 
with a loop $K$ based on $\ast \in M$. Let also $G$
be a degree $2$ clover in $M$ disjoint from $K$ and 
let $(M_G, K_G)$ be the result of the surgery along $G$.
Then, there exists an isomorphism
\begin{diagram}
\frac{\pi_1(M,\ast)}{\pi_1(M,\ast)_3} & \rTo^{\simeq} &
\frac{\pi_1(M_G,\ast)}{\pi_1(M_G,\ast)_3} 
\end{diagram}
sending $[K]$ to $[K_G]$.
\end{lemma}
\noindent
This lemma allows us to conclude with the following proposition-definition.
\begin{proposition}
For homology cylinders over $\Sigma=\Sigma_{g,1}$ or $\Sigma_g$, 
there are some well-defined homomorphisms
\begin{diagram}
\overline{\mathcal{C}}_{1}(\Sigma_{g,1})& \rTo^{\eta_1} & \Lambda^3H &
\textrm{ and } \ & 
\overline{\mathcal{C}}_{1}(\Sigma_{g})& \rTo^{\eta_1} & 
\frac{\Lambda^3H}{\omega \wedge H}.
\end{diagram}
Induced by the map (\ref{eq:brut}), they are called 
the \emph{first Johnson homomorphisms}.
\end{proposition}
\begin{remark} \label{rem:Johnson} 
The composition of $\eta_1$ with the map $C:\mathcal{T}(\Sigma)
\rTo \overline{\mathcal{C}}_1(\Sigma)$ is the classical  
homomorphism defined in \cite{Jab}.
\end{remark}
\noindent \textbf{Proof of Lemma \ref{lem:Y2}.}
Using \cite[Lem. 5.1]{GGP}, one shows that
$$M_G \cong_+ M\setminus int(N(G)) \ \cup_{j|_\partial} \ (H_4)_L$$
where $H_4\rInclus^j M$ is an oriented embedding of the standard
genus $4$ handlebody onto $N(G)$, which is a regular neighborhood of $G$
in $M$, and where $L$ is the $2$-component framed link 
shown\footnote{Blackboard framing convention is used.}\,\, 
on Fig. \ref{fig:link}. Through this diffeomorphism 
$K_G$ goes to  $K\subset M\setminus int(N(G))$.\\
Moreover, $L$ is Kirby-equivalent to the 
$3$-component link $N$ drawn on the right part 
of Fig. \ref{fig:link}.
It turns out that $N$ is a boundary link.
More precisely, up to a $(\pm 1)$-framing correction,
one can push disjointly $N_3$, $N_1$ and then $N_2$
to the boundary of $H_4$. We obtain some simple closed curves on 
$\Sigma_4=\partial H_4$, which are bounding curves.
Therefore, twist along each of these curves induces the identity at the level of 
$\pi_1(\Sigma_4,\ast)/\pi_1(\Sigma_4,\ast)_3$.
We then obtain the lemma by a Van-Kampen type argument.
$\square$\\

\subsection{Birman-Craggs homomorphisms for homology cylinders} 
Birman-Craggs homomorphisms were defined in \cite{BC} 
and they were enumerated 
in \cite{JBC}. Levine also outlined in \cite{L}
how they can be extended to homology cylinders.
In this paragraph, we review Birman-Craggs homomorphisms 
in a self-contained way. For this, we use the spin refinement of the Goussarov-Habiro 
theory of finite type invariants, introduced by the first author in \cite{M}.

We first fix a few notation. If $(M,\sigma)$ is a closed spin
$3$-manifold, let $R(M,\sigma)\in \mathbf{Z}_{16}$ denote its
Rochlin invariant. If $M$ is a homology sphere, we
will denote its (unique) spin structure by $\sigma_0$.
Recall from \cite{M} that surgery along a $Y$-graph makes also sense
among spin $3$-manifolds:
\begin{displaymath}
\left(
\begin{array}{cl}
\textsf{Data:} &  \textrm{(i)} \ (M,\sigma), \textrm{ a closed spin $3$-manifold} \\
&  \textrm{(ii)} \ G, \textrm{ a } Y \textrm{-graph in } M
\end{array} \right) 
\rightsquigarrow
\textsf{Result: } (M_{G},\sigma_{G}).
\end{displaymath}
The following lemma describes precisely how the Rochlin invariant
is modified during surgery along a $Y$-graph.
\begin{lemma}
Let $(M,\sigma)$ be a closed spin $3$-manifold, and
let $G$ be a $Y$-graph in $M$ whose leaves are ordered,
oriented and denoted  by $K_1$, $K_2$ and $K_3$. Then,
\begin{equation}
\label{eq:Rochlin}
R(M_G,\sigma_G)-R(M,\sigma)=8\cdot\prod_{k=1}^3 
e(t_{K_k})(\sigma)\in \mathbf{Z}_{16},
\end{equation}
where $8\cdot :\mathbf{Z}_2 \rInto \mathbf{Z}_{16}$ denotes
the usual injection, and where $t_{K_k} \in H_1(E(FM);\mathbf{Z})$
and $e(t_{K_k}) \in A(Spin(M),\mathbf{Z}_2)$ have been defined
in \S $2.2$.
\end{lemma}
\noindent \textbf{Proof.}
Let $j:H_3 \rInclus M$ be the embedding of the genus $3$ handlebody,
determined (up to isotopy) by the $Y$-graph $G$ in $M$. Then, it follows
from \cite[Prop. 1]{M}, that the variation 
$R(M_G,\sigma_G)-R(M,\sigma)$ only depends on $j^*(\sigma) \in Spin(H_3)$.
Also, according to equation (\ref{eq:Boolean}) from the proof 
of Lemma \ref{lem:type}, the rhs of (\ref{eq:Rochlin}) is determined
by $j^*(\sigma)\in Spin(H_3)$.\\ 
For $i_1,i_2,i_3 \in \{0,1\}$, 
we denote by $G_{i_1i_2i_3}$ the trivial $Y$-graph  
in $\mathbf{S}^3$ (with ordered and oriented leaves) and whose 
leaf number $k$ is trivial and $i_k$-framed; we also denote by 
$j_{i_1i_2i_3}: H_3 \rInclus \mathbf{S}^3$  the corresponding embedding. Then, 
\begin{displaymath}
Spin(H_3)=\left\{j_{i_1i_2i_3}^*(\sigma_0) \arrowvert 
i_1,i_2,i_3 \in \{0,1\}\right\}.
\end{displaymath}
Thus, it is enough to prove (\ref{eq:Rochlin}) 
when $(M,\sigma)$ is $\left(\mathbf{S}^3,\sigma_0\right)$ 
and when $G$ is a $G_{i_1i_2i_3}$, so that we now
restrict ourselves to this case.
By Lemma \ref{lem:type} b) (iv), the rhs of equation (\ref{eq:Rochlin}) is 
$8$ if $i_1=i_2=i_3=1$ and is $0$ otherwise.
The same holds for the lhs of equation (\ref{eq:Rochlin}).
Indeed, surgery along a $Y$-graph with a trivial leaf has no effect
(by the ``blow-up move'' of \cite[Fig. 6]{GGP}), and 
surgery  on $\mathbf{S}^3$ along $G_{111}$ gives the
Poincar\'e sphere whose Rochlin invariant 
is $8\in \mathbf{Z}_{16}$. It follows that equation (\ref{eq:Rochlin}) 
holds in these eight particular cases.
$\square$\\

Let $\Sigma$  be $\Sigma_g$ or $\Sigma_{g,1}$.
Let $j$ be an oriented embedding of $\Sigma$ in $\mathbf{S}^3$,
and let $M=(M,i^+,i^-)$ be a homology cylinder over $\Sigma$. We can then
cut $\mathbf{S}^3$ along $Im(j)$, and  glue back $M$ (using the 
identifications $j$, $i^+$ and $i^-$). We get a new homology sphere which is denoted by
\begin{displaymath}
\mathbf{S}^3(M,j).
\end{displaymath} 
It is shown in \cite[Cor. 1]{M} that the Rochlin invariant is a degree $1$
invariant: in particular, it is preserved under a $Y_2$-surgery.  
Therefore, $R\left(\mathbf{S}^3(M,j),\sigma_0\right)$ only depends on the 
$Y_2$-equivalence class of $M$ (and $j$). 
Suppose now we are given a surgery
presentation of the $Y_2$-equivalence class of $M$ on $1_\Sigma$:
\begin{displaymath}
\psi_1\left(\sum_{i=1}^n\Y \left[z_1^{(i)},z_2^{(i)},z_3^{(i)}\right]\right)
=M\in \overline{\mathcal{C}}_1(\Sigma).
\end{displaymath}
Recall that the labels $z_k^{(i)}$  belong to $P$
and thus give some $e\left(z_k^{(i)}\right)\in B_g^{(1)}$. We also put 
$\sigma=j^{*}(\sigma_0) \in Spin(\Sigma)$, 
which can be identified with the quadratic form
$q_\sigma \in \Omega_g$
according to the Johnson construction (see \S  2.2). 
We then deduce from (\ref{eq:Rochlin})
the following \emph{cubic} formula:
\begin{equation} \label{eq:cubic}
\frac{ R\left(\mathbf{S}^3(M,j),\sigma_0\right) }{ 8 } = \sum_{i=1}^n \prod_{k=1}^3
e\left(z_k^{(i)}\right)(q_\sigma) \in \mathbf{Z}_2.
\end{equation}
In particular, this shows that:
\begin{enumerate}
\item[(i)] $R\left(\mathbf{S}^3(M,j),\sigma_0\right)$ only depends on 
$\sigma=j^*(\sigma_0)\in Spin(\Sigma)$ 
(and the $Y_2$-equivalence class of $M$);
\item[(ii)] if $N$ is another homology cylinder over $\Sigma$, then:
\begin{displaymath}
\frac{R\left(\mathbf{S}^3(M\cdot N,j),\sigma_0\right)}{8}=
\frac{R\left(\mathbf{S}^3(M,j),\sigma_0\right)}{8}+
\frac{R\left(\mathbf{S}^3(N,j),\sigma_0\right)}{8} \in \mathbf{Z}_2.
\end{displaymath}
\end{enumerate}
We now distinguish the case $\Sigma=\Sigma_g$ from the case $\Sigma=\Sigma_{g,1}$.\\

\begin{figure}[!h] 
\begin{center} 
\includegraphics[width=8cm,height=4cm]{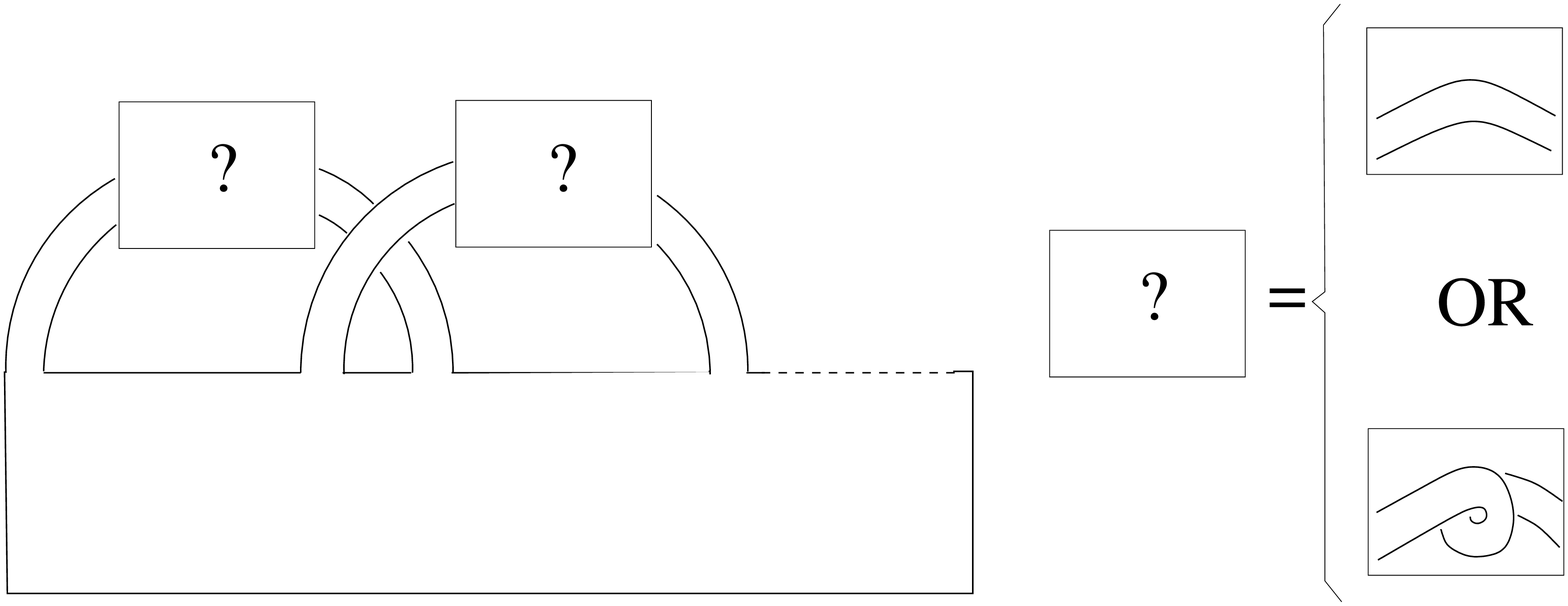} 
\fcaption{Some particular embeddings of $\Sigma_{g,1}$
in $\mathbf{S}^3$}
\label{fig:particular_embed_boundary} 
\end{center} 
\end{figure}

In the boundary case, any spin structure $\sigma$ on $\Sigma_{g,1}$ can
be realized as a $j^*(\sigma_0)$ for a certain embedding $j:\Sigma_{g,1}
\rInclus \mathbf{S}^3$. In fact, the specific embeddings of
$\Sigma_{g,1}$ whose images are depicted in Fig. 
\ref{fig:particular_embed_boundary} do suffice.\\
 
\begin{figure}[!h] 
\begin{center} 
\includegraphics[width=8cm,height=4cm]{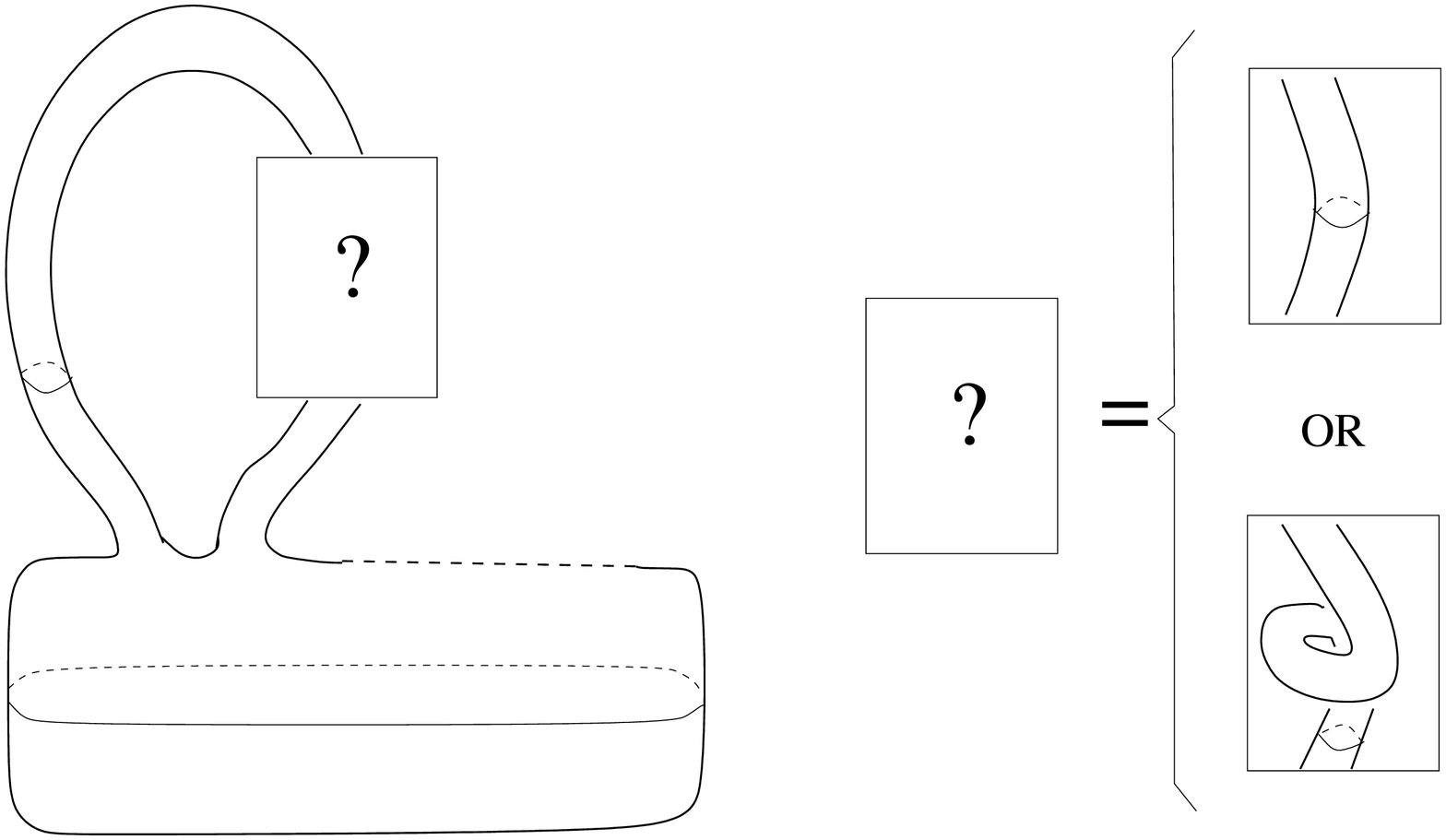}
\fcaption{Some particular embeddings of $\Sigma_{g}$
in $\mathbf{S}^3$}
\label{fig:particular_embed_closed} 
\end{center} 
\end{figure}

As for the closed case, observe that any
embedding $j:\Sigma_{g} \rInclus \mathbf{S}^3$ is splitting, 
so that $\sigma=j^*(\sigma_0)$ is spin-bounding. Conversely,
any spin structure on $\Sigma_g$ which spin-bounds can be
so realized: choose an appropriate embedding  of 
$\Sigma_g$ among the particular ones whose images are shown in Fig.
\ref{fig:particular_embed_closed}.
\\
Two other facts about these structures still have to be mentioned.
First, $\sigma\in Spin(\Sigma_g)$ spin-bounds if and only if
the Arf invariant $\alpha(q_\sigma)$ vanishes 
(see \cite[p.36]{Kirby}). Second, if $f$ and $f'$ are two 
cubic polynomials on $\Omega_g$ (namely $f,f' \in B^{(3)}_g$), then
they are identical on the quadratic forms with trivial Arf invariant
if and only if $f-f'$ is a multiple of $\alpha$ (see \cite[Lem. 14]{JBC} 
for a proof\footnote{There, the proof is given for a genus $g\ge 3$,
but the same arguments allow us to prove that this fact also holds 
for a genus $g=0,1$ or $2$.}\,\, of this algebraic fact).\\

All of our present discussion leads to the following proposition-definition.
\begin{proposition}
\label{prop:def_BC}
There exist some well-defined homomorphisms
\begin{diagram}
\overline{\mathcal{C}}_1(\Sigma_{g,1})& \rTo^\beta & B_g^{(3)} &
\textrm{ and }\ & 
\overline{\mathcal{C}}_1(\Sigma_{g})& \rTo^\beta & 
\frac{B_g^{(3)}}{\alpha \cdot B_g^{(1)}},
\end{diagram}
such that, for $M$ a homology cylinder over $\Sigma$
and for $j:\Sigma\rInclus \mathbf{S}^3$ an oriented embedding, we have
\begin{displaymath} 
\beta(M)\left(q_{j^*(\sigma_0)}\right)=
\frac{R\left(\mathbf{S}^3(M,j),\sigma_0\right)}{8} \in \mathbf{Z}_2.
\end{displaymath}
\end{proposition}
\begin{remark}
\label{rem:BC}
By composing $\beta$ with the map 
$C:\mathcal{T}(\Sigma)\rTo \overline{\mathcal{C}}_1(\Sigma)$, 
we obtain the classical Birman-Craggs homomorphisms, as 
presented by Johnson in \cite{JBC}.  
\end{remark} 

\section{Proof of the Results}
\label{sec:proof} 
In this section, we prove the results announced in the introduction. 
\begin{conventions}
In the proofs, we will use some specific techniques of Habiro.
Recall that its calculus of claspers developped in \cite{H}
is based on the definition of surgery along a \emph{basic clasper}. 
So as to be consistent with our Conventions \ref{conv:Y},
we define here  a \emph{basic clover} $C$ 
in a 3-manifold $M$ to be the embedding 
into $M$ of the surface depicted on the left part 
of Fig. \ref{fig:basic}.  
\emph{Surgery along $C$} is defined as the surgery
along the $2$-component framed link shown\footnote{Blackboard
framing convention is used.}\,\, in the right part 
of Fig. \ref{fig:basic}.
\begin{figure}[!h] 
\begin{center} 
\includegraphics{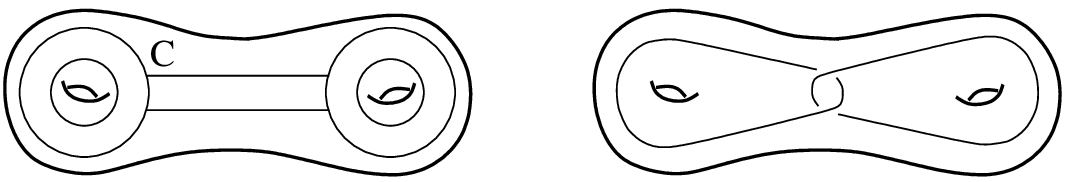} 
\fcaption{A basic clover $C$ 
and the associated framed link} 
\label{fig:basic} 
\end{center} 
\end{figure} 
Then, a basic clover is a basic clasper but with opposite 
surgery meaning. Consequently, before using one of the thirteen Habiro's moves,
\emph{we will have to take its mirror image}.
\end{conventions}

\subsection{$Y$-equivalence: proof of Proposition 
\ref{prop:homology_cylinders}}
\label{subsec:Y-equivalence}
Since surgeries along clovers preserve homology, the inclusions
$\mathcal{C}_1(\Sigma_g) \subset \mathcal{HC}(\Sigma_g)$ 
and $\mathcal{C}_1(\Sigma_{g,1}) \subset \mathcal{HC}(\Sigma_{g,1})$
are clear.\\

We now prove the inclusion $\mathcal{HC}(\Sigma_{g,1})\subset
\mathcal{C}_1(\Sigma_{g,1})$ using a result of Habegger. 
For this, we need the following definition. Let $k\geq 0$ be an integer, 
a \emph{homology handlebody} of genus $k$ is a pair $(M,i)$ where
\begin{itemize}
\item[(i)] $M$ is a compact oriented $3$-manifold whose integral
homology groups
are isomorphic to those of $H_k$, the standard genus $k$ handlebody;
\item[(ii)] $i: \Sigma_k=\partial H_k \rTo  M $ is an 
oriented embedding with image $\partial M$.
\end{itemize}
\begin{theorem}[Habegger, \cite{Ha}]
\label{th:handlebodies}
Let $(M_1,i_1)$ and $(M_2,i_2)$ be genus $k$ 
homology handlebodies such that
\begin{displaymath}
Ker\left(H_1\left(\Sigma_k;\mathbf{Z}\right)
\rTo^{i_{1,*}} H_1\left(M_1;\mathbf{Z}\right)\right)=
Ker\left(H_1\left(\Sigma_k;\mathbf{Z}\right)
\rTo^{i_{2,*}} H_1\left(M_2;\mathbf{Z}\right)\right).
\end{displaymath}
Then, $(M_1,i_1)$ and  $(M_2,i_2)$ are $Y$-equivalent.
\end{theorem} 
\noindent
In the sequel we identify $H_{2g}$ with $\Sigma_{g,1}\times I$,
and so $\Sigma_{2g}$ with $\partial\left(\Sigma_{g,1\times I}\right)$.
We also denote by $\left(H_{2g},j\right)$ the standard genus $2g$
handlebody, with inclusion $j:\Sigma_{2g}\rInclus H_{2g}$. 
Any homology cobordism $M=(M,i^+,i^-)$ over $\Sigma_{g,1}$
produces a genus $2g$ homology handlebody $(M,i)$, by defining
$i:\Sigma_{2g}\rTo M$ to be the diffeomorphism obtained from the gluing of
$i^+$ with $i^-$. Suppose now that $M$ is a homology cylinder.
Proving that the homology handlebody $(M,i)$ is $Y$-equivalent 
to $\left(H_{2g},j\right)$ will imply that the homology cylinder
$M$ is $Y$-equivalent to $\left(\Sigma_{g,1}\times I, Id,Id\right)$.\\
For this, let  $x_1^*,\ldots,x_g^*, y_1^*, \dots,y_g^*$ 
be some disjoint proper arcs in $\Sigma_{g,1}$, 
which are  ``dual'' to the loops $x_1,\ldots,x_g, y_1,\ldots,y_g$ of Fig. \ref{fig:basis}, 
in the sense that  $x_k^*$ (resp. $y_k^*$) transversely intersects $x_k$ (resp. $y_k$) 
once but does not intersect the other loops. 
For example, choose the first attaching region of each 1-handle.  
For each $k$, $X_k=x_k^*\times I$ and $Y_k=y_k^*\times I$ 
are discs in $\Sigma_{g,1}\times I$. 
The kernel of $j_*:H_1\left(\Sigma_{2g}\right)\rTo H_1\left(\Sigma_{g,1}\times I\right)$ is spanned by
$\partial X_1, \dots, \partial X_g, \partial Y_1,\dots, \partial Y_g$.
On the other hand, observe that $\pm \partial Y_k$ (resp. $\pm \partial X_k$) 
is homologous to $x_k \times 0- x_k\times 1$ 
(resp. to $y_k\times 0-y_k\times 1$) in $\Sigma_{2g}$. Therefore, since $M$ is a homology 
cylinder, $i\left(\partial X_k\right)$ and $i\left( \partial Y_k\right)$ are nul-homologous in $M$.
As the kernel of $i_*:H_1\left(\Sigma_{2g}\right) \rTo H_1(M)$ has to be of dimension $2g$,
it is spanned by  $\partial X_1, \dots, \partial X_g, \partial Y_1,\dots, \partial Y_g$.
It follows from Th. \ref{th:handlebodies} that $(M,i)$ is $Y$-equivalent to $\left(H_{2g},j\right)$,  
which proves the inclusion $\mathcal{HC}\left(\Sigma_{g,1}\right)\subset 
\mathcal{C}_1\left(\Sigma_{g,1}\right)$. \\

Let us now justify the inclusion $\mathcal{HC}\left(\Sigma_{g}\right)\subset
\mathcal{C}_1\left(\Sigma_{g}\right)$.
Let $j:\Sigma_{g,1}\rInclus \Sigma_g$ be an embedding and 
let $D\subset \Sigma_g$ be its complementary disk.
Take a homology cobordism $M=(M,i^+,i^-)$ over $\Sigma_{g,1}$.
Then, the embedding $(i^+)|_\partial \circ(j|_\partial)^{-1}
=(i^-)|_\partial \circ (j|_\partial)^{-1}: \partial D \rInclus \partial M$ 
can be stretched to an embedding $\partial D \times I \rInclus \partial M$.
The latter allows us to attach the $2$-handle $D\times I$ to $M$. This results
in a homology cylinder over $\Sigma_g$. 
We have thus defined a  \emph{filling-up} map
\begin{diagram}
\mathcal{C}\left(\Sigma_{g,1}\right) & \rTo^{j}  &
\mathcal{C}\left(\Sigma_{g}\right),
\end{diagram}
which is obviously surjective. 
Let $M\in \mathcal{HC}\left(\Sigma_{g}\right)$,
and pick a $N\in \mathcal{C}\left(\Sigma_{g,1}\right)$ 
such that $M$ is a filling-up
of $N$. Then, $N$ has to be a homology cylinder and so is $Y$-equivalent
to $1_{\Sigma_{g,1}}$. We conclude that 
$M\in\mathcal{C}_1\left(\Sigma_{g}\right)$,
which completes the proof of Proposition \ref{prop:homology_cylinders}. 

\subsection{The boundary case: proof of Theorem \ref{th:boundary}}
\label{subsec:bound}
Recall from Example \ref{ex:functors} that the Abelian group  
$\mathcal{A}_1(H,0)$ can be identified with $\Lambda^3H$, and  
likewise $\mathcal{A}_1\left(H_{(2)},0\right)$ with $\Lambda^3H_{(2)}$. 
The following lemma will allow us to identify
$\mathcal{A}_1\left(B_g^{(1)},\overline{1}\right)$ with $B_{g}^{(3)}$.
\begin{lemma} \label{lem:B1}
Let $\gamma: \mathcal{A}_1\left(B_g^{(1)},\ov{1}\right) \rTo   
B_{g}^{(3)}$ be the map given by 
multiplying the labels of the abstract
$Y$-graphs: $\gamma(\Y[z_1,z_2,z_3])=z_1z_2z_3$.
Then, $\gamma$ is a well-defined isomorphism. 
\end{lemma} 
\noindent \textbf{Proof.} 
The fact that $\gamma$ is well-defined is clear.
In order to show that $\gamma$ is an isomorphism, it suffices
to construct an epimorphism $B_{g}^{(3)} \rOnto^\epsilon 
\mathcal{A}_1\left(B_g^{(1)},\ov{1}\right)$ such that 
$\gamma\circ \epsilon$ is the identity.\\
By choosing a basis $(e_j)_{j=1}^{2g}$ for $H$, one
determines an isomorphism between 
$B_g^{(3)}$  and  $\mathbf{Z}_2\oplus H_{(2)} \oplus \Lambda^2 H_{(2)}
\oplus \Lambda^3 H_{(2)}$: for $k=1,2,3$
and $j_1,\dots, j_k \in\{1,\dots,2g\}$ pairwise distinct,
the monomial $\prod_{i=1}^k \overline{e_{j_i}}$ 
is identified with the wedge product
$\wedge_{i=1}^k e_{j_i}$, and $\overline 1$ with $1\in \mathbf{Z}_2$.
Since $B_g^{(1)}$ is a period $2$ group, so is   
$\mathcal{A}_1\left(B_g^{(1)},\ov{1}\right)$ 
by the multilinearity relation. 
Then, it suffices to define $\epsilon$ on the above mentioned 
$\mathbf{Z}_2$-basis of $\mathbf{Z}_2\oplus H_{(2)} \oplus \Lambda^2 H_{(2)}
\oplus \Lambda^3 H_{(2)} \simeq B_g^{(3)}$. We put
$\epsilon(1)=\Y\left[\overline{1},\overline{1},\overline{1}\right]$,
$\epsilon(e_j)=\Y\left[\overline{e_j},\overline{1},\overline{1}\right]$,
$\epsilon(e_{j_1}\wedge e_{j_2})=\Y\left[\ov{e_{j_1}},
\ov{e_{j_2}},\overline{1}\right]$ 
(with $j_1 \neq j_2$) and 
$\epsilon(e_{j_1}\wedge e_{j_2} \wedge e_{j_3})=
\Y\left[\ov{e_{j_1}},\ov{e_{j_2}},\ov{e_{j_3}}\right]$ 
(with $j_1,j_2,j_3$ pairwise distinct). The map
$\epsilon$ is surjective by the multilinearity 
and slide relations, and obviously
satisfies $\gamma\circ \epsilon=Id$. 
$\square$\\

Recall from \S 2.2 that the maps
\begin{diagram}
P & \rTo^p & (H,0) & \textrm{ and } & 
P & \rTo^e & \left(B_g^{(1)},\overline{1}\right)
\end{diagram}
are the canonical projections of the pullback of special Abelian groups
\begin{displaymath}
P=(H,0)\times_{(H_{(2)},0)}\left(B_g^{(1)},\ov{1}\right).
\end{displaymath}
They happen to be surjective.
\begin{lemma} \label{lem:john}
The following diagram commutes: 
\begin{diagram}
\mathcal{A}_1(P) & \rOnto^{\psi_1} & \ov{\mathcal{C}}_1(\Sigma_{g,1}) \\
& \rdOnto<{\mathcal{A}_1(p)} & \dTo>{\eta_1} \\
& & \mathcal{A}_1(H,0).
\end{diagram} 
\end{lemma}

\noindent \textbf{Proof.}
Let us verify that $\eta_1\left(\psi_1(Y)\right)
=\mathcal{A}_1\left(p\right)(Y)$ for a generator
$Y=\Y[z_1,z_2,z_3]$ of $\mathcal{A}_1(P)$.
We put $M=\psi_1(Y)$,  so that $M=(1_{\Sigma_{g,1}})_G$
where $G$ is an appropriate $Y$-graph as 
described in \S 2.3. Its leaves are in particular ordered and oriented,
they are denoted by $K_1$, $K_2$ and $K_3$: $[K_i]=p(z_i) \in H$.
Set $\pi=\pi_1(\Sigma_{g,1},\ast)$ and 
let $\overline y \in \pi/\pi_3$ be represented by $y\in \pi$:
we want to compute $\eta_1\left(M\right)$ on $\overline{y}$. 
This goes as follows: 
choose an immersed based curve $k$ in $\Sigma_{g,1}^+$
representing $y$ (\emph{via} the identification of
$\Sigma_{g,1}$ with $\Sigma_{g,1}^+$),  pick an 
oriented based knot $K \subset M $ in a collar
of $\Sigma_{g,1}^+$ which is a push-off of $k$, and  find another based knot 
$K'\subset M$ in a collar of $\Sigma_{g,1}^-$ such that
the pairs $(M,K)$ and $(M,K')$ are $Y_2$-equivalent.
Then (\emph{via} the identification of
$\Sigma_{g,1}$ with $\Sigma_{g,1}^-$), this knot $K'$ determines a $y'\in \pi$,
and by Lemma \ref{lem:Y2}, 
the result $\eta_1(M)(\overline{y})$ 
is then $\overline{y'}\in \pi/\pi_3$.
We now explain the procedure how to construct $K'$ from $K$.
 
\begin{figure}[!h] 
\begin{center} 
\includegraphics{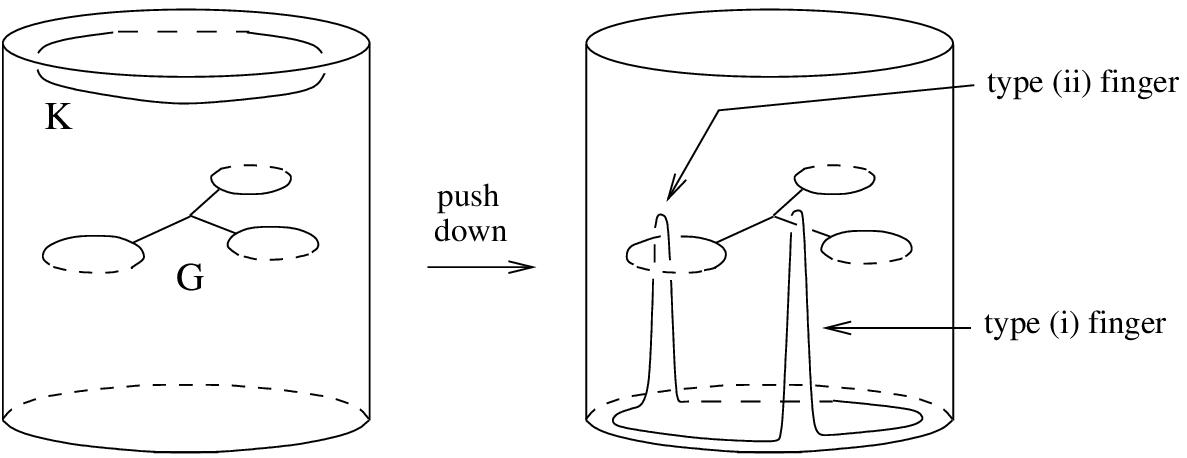}
\fcaption{Pushing the curve $K$ down the cylinder} \label{fig:fing} 
\end{center} 
\end{figure}

In $1_{\Sigma_{g,1}}\setminus G$, $K$ can be pushed down in a collar
of $\Sigma_{g,1}^-$ \emph{up to} some ``fingers'' which are of two types (see Fig. \ref{fig:fing}):
\begin{enumerate}
\item[(i)] the finger is pointing on an edge of $G$,
\item[(ii)] the finger is pointing on an leaf $K_i$ of $G$.
\end{enumerate}
But, each finger of type (i) can be isotoped along the corresponding edge
towards its leaf and so can be replaced by two fingers of type (ii), so that 
up to some isotopy of the immersed curve $k$ in $\Sigma_{g,1}^+$,
we can suppose each finger to be of type (ii).
Since $K_i$ has been oriented, each finger comes with a sign. 
Let $k_i$ be an immersed curve on $\Sigma_{g,1}^+ \subset 1_{\Sigma_{g,1}}$
such that $[k_i]=p(z_i) \in H$. We can suppose that $K_i$ is 
a push-off of $k_i$ (with possibly an additional twist):
there are then as many fingers 
as intersection points of $k_i$ with $k$ in $\Sigma_{g,1}^+$;
the sign of the finger corresponds with the sign of the intersection
point contributing to $[k]\bullet [k_i] \in \mathbf{Z}$.\\ 
A finger move can be realized by surgery on a basic clover.
Let $K'$ be a copy of $K$ in a
collar of $\Sigma_{g,1}^-\subset 1_{\Sigma_{g,1}}\setminus G$.
There is then a family of basic clovers 
$\left(C^{(i)}_j\right)_{j=1,\dots,n_i}^{i=1,2,3}$ 
in $1_{\Sigma_{g,1}}\setminus G$, 
such that each $C_j^{(i)}$ has a simple leaf which laces $K'$
and another simple leaf wich laces the leaf $K_i$, and such that:
\begin{displaymath}
(M,K) \textrm{ is diffeomorphic to } 
\left( 1_{\Sigma_{g,1}},K' 
\right)_{\left(\cup_{i,j} C_j^{(i)}\right) \cup G}.
\end{displaymath}
According to the sign of the corresponding finger, each basic clover comes 
with a sign denoted by $\varepsilon(i,j)$.
\emph{Cutting} the leaf $K_1$ (see \cite[Cor. 4.3]{GGP})
$n_1$ times, we obtain $n_1$ new $Y$-graphs $G_j^{(1)}$ 
($j \in \{1,\dots,n_1\}$): two leaves of $G_j^{(1)}$ are copies of $K_2$  
and $K_3$, and the third leaf forms with a leaf of $C_j^{(1)}$ the 
Hopf-link. Hence, by applying Habiro move $2$ (or \cite[Th. 2.4]{GGP}) 
to $C_j^{(i)}\cup G_j^{(i)}$ we obtain a new $Y$-graph still denoted 
by $G_j^{(i)}$. We do the same for $i=2$ and $i=3$,
therefore:
\begin{displaymath}
(M,K) \textrm{ is } Y_2\textrm{-equivalent to } 
\left( 1_{\Sigma_{g,1}},K' 
\right)_{\left(\cup_{i,j}G_j^{(i)}\right)\cup G}.                       
\end{displaymath}
Up to $Y_2$-equivalence of the pair 
$\left( 1_{\Sigma_{g,1}},K'\right)_{G\cup \left(\cup_{i,j} G_j^{(i)}\right)}$,
one can suppose that, for each $(i,j)$, the whole of $G_j^{(i)}$
lies in a collar neighborhood of 
$\Sigma_{g,1}^- \subset 1_{\Sigma_{g,1}}$. We now do the 
surgery along $G$, and then along each of the $G_j^{(i)}$:
the latter does not modify the
$3$-manifold $M$ but changes the knot. The new knot we obtain is
still denoted by $K'$ and satisfies the announced required properties.

We now calculate the $y' \in \pi$ defined by $K'$.
In view of Habiro move 10, the contribution of each $Y$-graph
$G_j^{(1)}$ to the modification of $K'$ is in $\pi$
the commutator $\left[k_2,k_3^{-1}\right]^{\varepsilon(1,j)}$. Therefore, we obtain
\begin{equation}
y'\cdot y^{-1}= \prod_{i\in \mathbf{Z}_3}
\left[k_{i+1},k_{i+2}^{-1}\right]^{[k]\bullet [k_i]}
\in \frac{\pi_2}{\pi_3}.
\end{equation}
Then, as a homomorphism $H\rTo \pi_2/\pi_3=L_2(H)$, $\eta_1(M)$
sends any $h\in H$ to
\begin{displaymath}
-\sum_{i\in \mathbf{Z}_3} \left(h\bullet p(z_i)\right)
\cdot [p(z_{i+1}),p(z_{i+2})] \in L_2(H).
\end{displaymath}
which corresponds to $\sum_{i\in \mathbf{Z}_3}p(z_i)\otimes 
[p(z_{i+1}, p(z_{i+2})]$ in $H\otimes L_2(H)$, 
to $p(z_1)\wedge p(z_2)\wedge p(z_3)$ in $\Lambda^3H$, and 
so to $\mathcal{A}_1(p)(Y)$.
$\square$\\

\begin{lemma} \label{lem:BC}
The following diagram commutes:
\begin{diagram} 
\mathcal{A}_1(P) & \rOnto^{\psi_1} & \ov{\mathcal{C}}_1(\Sigma_{g,1}) \\  
 & \rdOnto<{\mathcal{A}_1(e)} & \dTo>{\beta} \\ 
& & \mathcal{A}_1\left(B_g^{(1)},\ov{1}\right).
\end{diagram} 
\end{lemma}
\noindent \textbf{Proof.}
According to the definition of $\beta$ we gave in
Proposition \ref{prop:def_BC}, this is a direct consequence
of equation (\ref{eq:cubic}).
$\square$\\

We still denote by
\begin{diagram}
(H,0) & \rTo^{-\otimes\mathbf{Z}_2} & \left(H_{(2)},0\right) &
\textrm{ and } & \left(B^{(1)}_g,\overline{1}\right) 
& \rTo^{\kappa} & \left(H_{(2)},0\right), 
\end{diagram}
the maps which appear in the pullback diagram for $P$ 
(see \S 2.2). Then, as a consequence 
of the two preceeding lemmas, 
$\mathcal{A}_1(\kappa)\beta \psi_1 = \mathcal{A}_1(\kappa e) =
\mathcal{A}_1\left((-\otimes\mathbf{Z}_2)p \right)=
\mathcal{A}_1(-\otimes\mathbf{Z}_2)\eta_1 \psi_1$.
Since $\psi_1$ is an epimorphism, we get: 
$\mathcal{A}_1(\kappa)\beta = 
\mathcal{A}_1(-\otimes\mathbf{Z}_2)\eta_1$.
Construct the following pull-back: 
\begin{diagram}
\mathcal{A}_1(H,0) \times_{\mathcal{A}_1(H_{(2)},0)}\SEpbk 
\mathcal{A}_1\left(B_g^{(1)},\overline{1}\right)
& \rTo &  \mathcal{A}_1\left(B_g^{(1)},\overline{1}\right)\\ 
 \dTo &  & \dTo>{\mathcal{A}_1(\kappa)} \\ 
\mathcal{A}_1(H,0) & \rTo_{\mathcal{A}_1(-\otimes\mathbf{Z}_2)} 
& \mathcal{A}_1(H_{(2)},0)  
\end{diagram} 
which, through the above mentioned identifications, is
essentially the pull-back diagram for
$\Lambda^{3} H \times_{\Lambda^{3} H_{(2)}} B^{(3)}_g$
appearing in \S 1.3.
By the universal property of the pull-backs, there is then a homomorphism 
\begin{diagram}
\ov{\mathcal{C}}_1(\Sigma_{g,1})& \rTo^{(\eta_1,\beta)}& \mathcal{A}_1(H,0)
\times_{\mathcal{A}_1(H_{(2)},0)} 
\mathcal{A}_1\left(B_g^{(1)},\overline{1}\right)
\simeq  \Lambda^{3} H \times_{\Lambda^{3} H_{(2)}} B^{(3)}_g. 
\end{diagram}
Moreover, we also have by functoriality another natural map  
\begin{diagram} 
\mathcal{A}_1 \big(\underbrace{(H,0)\times_{(H_{(2)},0)}  
(B_g^{(1)},\overline{1})}_{P}\big) & \rTo^{\rho} & \mathcal{A}_1(H,0) 
\times_{\mathcal{A}_1(H_{(2)},0)} 
\mathcal{A}_1\left(B_g^{(1)},\overline{1}\right).
\end{diagram} 
Lemma \ref{lem:john} and Lemma \ref{lem:BC} can then be summarized
in the commutativity of the following diagram: 
\begin{diagram} 
\mathcal{A}_1(P) & \rOnto^{\psi_1} & \ov{\mathcal{C}}_1(\Sigma_{g,1}) \\  
& \rdTo<{\rho} & \dTo>{(\eta_1,\beta)} \\ 
& & \Lambda^{3} H \times_{\Lambda^{3} H_{(2)}} B^{(3)}_g. 
\end{diagram} 
The following lemma will be the final step in proving
Theorem \ref{th:boundary}.
\begin{lemma} 
\label{lem:rho}
The map $\rho: \mathcal{A}_1(P) \rTo \Lambda^{3} H 
\times_{\Lambda^{3} H_{(2)}} B^{(3)}_g$ 
is an isomorphism. 
\end{lemma} 
Assume Lemma \ref{lem:rho}. Then, from the previous commutative diagram, 
it follows that $\psi_1$ is injective, and so is an isomorphism: as a consequence,
the same holds for $(\eta_1,\beta)$. The commutativity of
\begin{diagram}
\overline{\mathcal{C}}_1(\Sigma_{g,1}) & \lTo^{C} & 
\frac{\mathcal{T}_{g,1}}{\mathcal{T}_{g,1}'}\\ 
\dTo<{(\eta_{1},\beta)}>{\simeq} & 
\ldTo>{(\eta_{1},\beta)} & \\ 
\Lambda^{3} H \times_{\Lambda^{3} H_{(2)}} B^{(3)}_g & & 
\end{diagram}
follows from Remark \ref{rem:Johnson} and Remark \ref{rem:BC}.
In particular, when $g\geq 3$, $C$ is an isomorphism because
$(\eta_1,\beta): \mathcal{T}_{g,1}/\mathcal{T}_{g,1}'\rTo 
\Lambda^{3} H \times_{\Lambda^{3} H_{(2)}} B^{(3)}_g $ is so by \cite{JFinal}.\\

\noindent \textbf{Proof of Lemma \ref{lem:rho}.}  
We proceed as in Lemma \ref{lem:B1}. It suffices to construct an epimorphism
\begin{diagram}
\Lambda^{3} H \times_{\Lambda^{3} H_{(2)}} B^{(3)}_g  &
\rOnto^{\epsilon} &  \mathcal{A}_1(P)
\end{diagram}
such that $\rho \circ \epsilon$ is the identity.\\
Pick a basis $(e_i)_{i=1}^{2g}$ of $H$: we have seen in the proof 
of Lemma \ref{lem:B1} that this choice determines an isomorphism between
$B_g^{(3)}$ and $\Lambda^3 H_{(2)}\oplus \Lambda^2H_{(2)}
\oplus H_{(2)}\oplus \mathbf{Z}_2$.
Thus, it also defines an isomorphism between  
$ \Lambda^{3} H \times_{\Lambda^{3} H_{(2)}} B^{(3)}_g $
and  $\Lambda^3H\oplus \Lambda^2H_{(2)}\oplus H_{(2)}\oplus \mathbf{Z}_2$.
We now define $\epsilon$ by putting
\begin{enumerate}
\item[(i)] $\epsilon(e_i\we e_j\we e_k)=
\Y\left[(e_i,\ov{e_i}),(e_j,\ov{e_j}),(e_k,\ov{e_k})\right]$,
with $1\leq i < j < k \leq 2g$,
\item[(ii)] $\epsilon(e_i\we e_j)=\Y\left[(e_i,\ov{e_i}),(e_j,\ov{e_j}),
(0,\ov{1})\right]$, with $1\leq i < j \leq 2g$,
\item[(iii)] $\epsilon(e_i)=\Y\left[(e_i,\ov{e_i}),(0,\ov{1}),(0,\ov{1})\right]$,
with $1\leq i \leq 2g$, 
\item[(iv)] and $\epsilon(1)=\Y\left[(0,\ov{1}),(0,\ov{1}),(0,\ov{1})\right]$.
\end{enumerate}
Here, elements of $P$ are denoted as in Remark \ref{rem:elts_of_P}.
This assignation well defines $\epsilon$ because (i) determines
$\epsilon$ on a basis of the free group $\Lambda^{3} H$,
while (ii),(iii) and (iv) assign elements of $\mathcal{A}_1(P)$ 
of order at most $2$ to each element basis of the 
$\mathbf{Z}_2$-vector space 
$\Lambda^2H_{(2)}\oplus H_{(2)}\oplus \mathbf{Z}_2$.
Obviously, $\epsilon$ followed by $\rho$ gives the identity. 
Take now any generator $\Y[z_1,z_2,z_3]$ of $\mathcal{A}_1(P)$. 
For $i=1,2,3$, $z_i\in P$ can be written as a linear combination
of some $(e_j,\overline{e_j})$ and $(0,\ov{1})$.
The multilinearity, AS and slide relation allow us to conclude
that $\Y[z_1,z_2,z_3]$ is realized by $\epsilon$. 
Thus, $\epsilon$ is surjective.
$\square$\\

\subsection{The closed case: proof of Theorem \ref{th:closed}}
An isomorphism
\begin{diagram}
\mathcal{A}_1(P)& \rTo^\rho &
\Lambda^{3} H \times_{\Lambda^{3} H_{(2)}} B^{(3)}_g
\end{diagram}
is defined formally in the same way as in the boundary case
(see Lemma \ref{lem:rho}).
Recall that S stands for the subgroup of the 
pullback $\Lambda^{3} H \times_{\Lambda^{3} H_{(2)}} B^{(3)}_g$ 
corresponding to $\omega\wedge H \subset \Lambda^{3} H$ and 
$\alpha\cdot B^{(1)}_g \subset B^{(3)}_g$.
Then, $\rho^{-1}(\textrm{S})$ is the subgroup of $\mathcal{A}_1(P)$ 
comprising the elements
\begin{displaymath}
\sum_{i=1} ^g\Y\left[(x_i,\overline{x_i}),(y_i,\overline{y_i}),z\right], \quad
\textrm{where } z \textrm{ is any element of } P. 
\end{displaymath} 
\begin{lemma} \label{lem:S} 
In the closed case, the surgery map $\psi_1$ defined in 
\S $2.3$ vanishes on the subspace $\rho^{-1}(\textrm{S})$. 
\end{lemma} 
\noindent
As mentioned in the introduction, these symplectic relations 
$\rho^{-1}(\textrm{S})$ appears in \cite{H} for higher degrees.\\

\begin{figure}[!h] 
\begin{center} 
\includegraphics{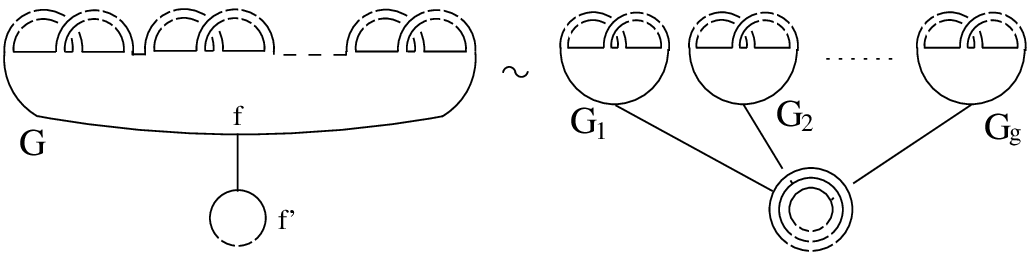}
\fcaption{Splitting the null-homologous leaf $f$} 
\label{fig:strut} 
\end{center} 
\end{figure}

\noindent \textbf{Proof of Lemma \ref{lem:S}.} 
Let $z \in P$, we aim to show that
\begin{equation}
\label{eq:to_be_proved}
\sum_{i=1}^g\psi_1\left(\Y\left[(x_i,\overline{x_i}),
(y_i,\overline{y_i}),z\right]\right) = 0 \in 
\overline{\mathcal{C}}_1(\Sigma_g).
\end{equation}
Consider in $1_{\Sigma_g}$ a basic clover $G$ with one trivial leaf $f$, 
and the other leaf $f'$ satisfying $t_{f'}=z\in P$.
Then, $f$ being trivial,
$\left(1_{\Sigma_g}\right)_G$ is diffeomorphic to $1_{\Sigma_g}$.
Furthermore, $f$ can be seen as a push-off of $\partial D$
where $D$ is a $2$-disk in $\Sigma_g^+$: in particular,
$f$ bounds the push-off of $\Sigma_g^+ \setminus D$ which is an 
embedded genus $g$ surface.
By applying Habiro moves 7 and 5, $f$ can be 
split in $g$  pieces so that $G$ is equivalent 
to the union of $g$ basic clovers denoted by $G_1,\dots, G_g$.
See Fig. \ref{fig:strut}. 
Each clover $G_i$ has a leaf which bounds a genus $1$ surface; 
by applying Habiro's move 10, it is  seen to be equivalent
to a $Y$-graph $G'_i$. According to Remark \ref{rem:elts_of_P},
the leaves of $G'_i$ represent
$(x_i,\overline{x_i})$, $(y_i,\overline{y_i})$ and $z$ in $P$,
so that $\left(1_{\Sigma_g}\right)_{G'_i}=
\psi_1\left(\Y\left[(x_i,\overline{x_i}),
(y_i,\overline{y_i}),z\right]\right)\in \overline{\mathcal{C}}_1(\Sigma_g)$.
Equation (\ref{eq:to_be_proved}) then follows.
$\square$\\

By the same arguments, appropriate versions of 
Lemma \ref{lem:john} and Lemma \ref{lem:BC} 
hold in the closed case: $\mathcal{A}_1(p)=\eta_1\circ 
\psi_1$ and $\mathcal{A}_1(e)=\beta\circ\psi_1$.
This leads us to a commutative diagram
\begin{diagram}  
 \frac{\mathcal{A}_1(P)}{\rho^{-1}(\textrm{S})} \quad & \rOnto^{\psi_1}   
&\overline{\mathcal{C}}_1(\Sigma_{g}) & \\ 
 & \rdTo<{\rho}>{\simeq}  & \dTo>{(\eta_{1},\beta)}& \\ 
& & \quad \frac{\Lambda^{3} H \times_{\Lambda^{3} H_{(2)}} 
B^{(3)}_g}{\textrm{S}} &\ \simeq\  
\frac{\Lambda^3 H}{\omega \wedge H} 
\times_{\left(\frac{\Lambda^3 H_{(2)}}{\omega_{(2)}\wedge H_{(2)}}\right)} 
\frac{B_g^{(3)}}{\alpha\cdot B_g^{(1)}}
\end{diagram} 
from which it follows that $\psi_1$, and then $(\eta_1,\beta)$, 
are isomorphisms. The commutativity of the right triangle
in Th. \ref{th:closed} is still given by Rem. \ref{rem:Johnson} and Rem. \ref{rem:BC}.\\

\subsection{Finite type invariants of degree $1$: proof of
Corollary \ref{cor:jbc}.} 
The equivalence (a)$\Leftrightarrow$(b) immediately results from the existence 
of the universal degree one additive invariant $\upsilon$ introduced 
in Remark \ref{rem:facts}.
The equivalence (c)$\Leftrightarrow$(a) is a direct consequence of
Theorem \ref{th:closed} and Theorem \ref{th:boundary}. 

\subsection{From the boundary case to the closed case}
In this last paragraph, we \emph{fix} an isomorphism 
\begin{diagram}
H_1(\Sigma_{g,1};\mathbf{Z}) & \rTo^\phi & H_1(\Sigma_{g};\mathbf{Z}).
\end{diagram}
It allows us to identify the sets 
$H$, $\Omega_g \simeq Spin(\Sigma)$, $B_g$ and $P$ corresponding
to $\Sigma_{g,1}$ with those of $\Sigma_{g}$.\\
Moreover, let $j:\Sigma_{g,1}\rInclus \Sigma_g$ be an embedding 
such that $j_*=\phi$ at the level of $H_1(-;\mathbf{Z})$. 
Recall from \S 4.1 the filling-up map, which can be restricted to
\begin{diagram}
\mathcal{C}_1\left(\Sigma_{g,1}\right) & \rTo^j & 
\mathcal{C}_1\left(\Sigma_{g}\right).
\end{diagram}
Note that it is compatible with the ``extending
by the identity'' map $\mathcal{T}_{g,1} \rTo \mathcal{T}_g$ defined by $j$,
and that it induces a group homomorphism 
$\overline{\mathcal{C}}_1\left(\Sigma_{g,1}\right)  \rTo  
\overline{\mathcal{C}}_1\left(\Sigma_{g}\right)$.
The latter can be verified to be independent on the choice of 
the embedding $j$ such that $j_*=\phi$, and so can be denoted by 
\begin{diagram}
\overline{\mathcal{C}}_1\left(\Sigma_{g,1}\right) & \rTo^\phi &
\overline{\mathcal{C}}_1\left(\Sigma_{g}\right).
\end{diagram}
The commutativity of the following diagram is easily proved from the 
various definitions:\\

\begin{diagram}
\mathcal{A}_1(P) & &\rTo^{\psi_1}  &  &
\overline{\mathcal{C}}_1(\Sigma_{g,1})  & & \\
&\rdTo_\rho & &\ldTo<{(\eta_{1},\beta)}  & \dNoto &\luTo>C & \\
\dOnto & & \Lambda^{3} H \times_{\Lambda^{3} H_{(2)}} B^{(3)}_g &
\lTo_{(\eta_{1},\beta)} &\HonV&  &\frac{\mathcal{T}_{g,1}}{\mathcal{T}_{g,1}'}\\
& &\dOnto & & \dOnto>\phi& & \\
\frac{\mathcal{A}_1(P)}{\rho^{-1}(\textrm{S})}\quad & 
\rNoto& \VonH &\rTo^{\psi_1} &\overline{\mathcal{C}}_1(\Sigma_{g}) & &\dOnto\\
& \rdTo_\rho& &\ldTo<{(\eta_{1},\beta)} & &\luTo^C & \\
& & \quad \frac{\Lambda^{3} H \times_{\Lambda^{3} H_{(2)}}  B^{(3)}_g}
{\textrm{S}}\quad &\lTo_{(\eta_{1},\beta)} & & &
\frac{\mathcal{T}_{g}}{\mathcal{T}_g'}\\ 
\end{diagram} 
\\
\\
\textbf{Acknowledgments} 
\\

The authors are indebted to Christian Blanchet and 
Nathan Habegger for helpful comments. Thanks are also
due to the anonymous referee for valuable suggestions.
\\
\\
\textbf{References} 
\\
 
\vspace{0.5cm}
\noindent
\footnotesize{Commutative diagrams were drawn with Paul Taylor's package.}
\end{document}